\documentclass{amsart}
\usepackage{amsmath}
\usepackage{amsxtra}
\usepackage{amscd}
\usepackage{amsthm}
\usepackage{amsfonts}
\usepackage{amssymb}
\usepackage{eucal}

\newcommand{\wt}{{\rm wt}\,} 
\newcommand{\gr}{{\rm gr}\,} 
\newcommand{\ch}{{\rm ch}}
\newcommand{\Xc}{\mathcal{X}}
\newcommand{\X}{\mathcal{X}}
\newcommand{\U}{U_{\sqrt{-1}}}
\newcommand{\Ures}{U_q^{{\rm res}}}
\newcommand{\qbin}[2]{{\left[
\begin{matrix}{\displaystyle #1}\\
{\displaystyle #2}\end{matrix}
\right]
}}

\newcommand{\nn}{\nonumber}
\newcommand{\bea}{\begin{eqnarray}}
\newcommand{\ena}{\end{eqnarray}}
\newcommand{\be}{\begin{eqnarray*}}
\newcommand{\en}{\end{eqnarray*}}

\newcommand{\C}{{\mathbb C}}
\newcommand{\Z}{{\mathbb Z}} 

\newcommand{\sltt}{\widetilde{\mathfrak{sl}}_2}
\newcommand{\res}{{\rm res}}
\newcommand{\Ch}{\widehat{C}}
\newcommand{\Rh}{\widehat{R}}
\newcommand{\Wh}{\widehat{W}}
\newcommand{\Dh}{\widehat{D}}
\newcommand{\Zc}{\mathcal{Z}}
\newcommand{\Zhc}{\widehat{\mathcal{Z}}}
\newcommand{\Zbc}{\overline{\mathcal{Z}}}
\newcommand{\Lc}{\mathcal{L}}
\newcommand{\Ps}{P^*}

\numberwithin{equation}{section}
\newtheorem{thm}{Theorem}[section]
\newtheorem{prop}[thm]{Proposition}
\newtheorem{lem}[thm]{Lemma}
\newtheorem{cor}[thm]{Corollary}

\begin{document} 

\title[]
{Form factors and action of
$U_{\sqrt{-1}}(\widetilde{\mathfrak{sl}}_2)$ on $\infty$-cycles}
\author{M. Jimbo, T. Miwa, E. Mukhin and Y. Takeyama}
\address{MJ: Graduate School of Mathematical Sciences, The
University of Tokyo, Tokyo 153-8914, Japan}\email{jimbomic@ms.u-tokyo.ac.jp}
\address{TM: Department of Mathematics, Graduate School of Science,
Kyoto University, Kyoto 606-8502
Japan}\email{tetsuji@math.kyoto-u.ac.jp}
\address{EM: Department of Mathematics,
Indiana University-Purdue University-Indianapolis,
402 N.Blackford St., LD 270,
Indianapolis, IN 46202}\email{mukhin@math.iupui.edu}
\address{YT: Department of Mathematics, Graduate School of Science,
Kyoto University, Kyoto 606-8502
Japan}\email{takeyama@kusm.kyoto-u.ac.jp}

\begin{abstract}
Let 
${\bf p}=\{P_{n,l}\}_{n,l\in\Z_{\ge 0}\atop n-2l=m}$ 
be a sequence of 
skew-symmetric polynomials in $X_1,\cdots,X_l$
satisfying $\deg_{X_j}P_{n,l}\le n-1$,  
whose coefficients are 
symmetric Laurent polynomials in $z_1,\cdots,z_n$.   
We call ${\bf p}$ an $\infty$-cycle if \par
\noindent$P_{n+2,l+1}\bigl|_{X_{l+1}=z^{-1},z_{n-1}=z,z_n=-z}$ 
$=z^{-n-1}\prod_{a=1}^l(1-X_a^2z^2)\cdot P_{n,l}$ 
holds for all $n,l$.   

These objects arise in integral representations for 
form factors of massive integrable field theory,
i.e., the $SU(2)$-invariant Thirring model and
the sine-Gordon model.
The variables $\alpha_a=-\log X_a$ are the integration variables
and $\beta_j=\log z_j$ are the rapidity variables.
To each $\infty$-cycle there corresponds a form factor
of the above models.  
Conjecturally all form-factors
are obtained from the $\infty$-cycles.


In this paper, we define an 
action of $U_{\sqrt{-1}}(\widetilde{\mathfrak{sl}}_2)$ 
on the space of $\infty$-cycles. 
There are two sectors of $\infty$-cycles depending on 
whether $n$ is even or odd. 
Using this action, we show 
that the character of the space of even 
(resp. odd) $\infty$-cycles 
which are polynomials in $z_1,\cdots,z_n$
is equal to the level $(-1)$ 
irreducible character of  
$\widehat{\mathfrak{sl}}_2$ 
with lowest weight $-\Lambda_0$ 
(resp. $-\Lambda_1$).
We also suggest a possible tensor product structure 
of the full space of $\infty$-cycles.
\end{abstract}

\maketitle

\setcounter{section}{0}
\setcounter{equation}{0}

\section{Introduction}

First let us recall the form factor bootstrap approach
to massive integrable models \cite{S1}
in field theory.
A form factor is a tower of meromorphic functions
\[
{\bf f}=(f_n(\beta_1,\ldots,\beta_n))_{n\geq0}
\]
in the variables $\beta_1,\ldots,\beta_n\in\C$,
satisfyng a certain set of axioms, 
to be referred to as Axiom $1$--$3$.
There are two sectors in the space of form factors:
the even sector where $f_n=0$ for all odd $n$, 
and the odd sector where $f_n=0$ for all even $n$.  
With each form factor is associated
a local field in the theory. Axiom 1 describes
the exchange relation of $\beta_j$ and $\beta_{j+1}$ in $f_n$, 
and Axiom 2 relates the analytic continuation 
$f_n(\beta_1,\ldots,\beta_{n-1},\beta_n+2\pi i)$
to the cyclic shift $f_n(\beta_n,\beta_1,\ldots,\beta_{n-1})$. 
These two axioms imply that, for each $n$, 
$f_n(\beta_1,\ldots,\beta_n)$ is a solution to the  
quantum Knizhnik-Zamolodchikov (qKZ) equation. 
Axiom 3 stipulates that $f_{n-2}$ is determined by the residue of $f_n$
at the simple pole $\beta_n=\beta_{n-1}+\pi i$.
In this paper we consider 
the $SU(2)$ invariant Thirring model (ITM). 
For the details of the axioms in this case, 
see Section \ref{formfactor}. 

One of the basic issues
in the theory of form factors is to 
describe all form factors satisfying the three axioms.
In recent papers \cite{N,N2}, Nakayashiki solved this problem
for the case of ITM under some assumptions. 
Subsequently, in \cite{JMT}
the case of the restricted sine-Gordon model 
was studied by a different method
based on representation theory of quantum affine algebras.
For the purpose of introducing the subject, 
and for making a comparison between the two methods, 
let us briefly review the results of \cite{N,N2,JMT}.

Nakayashiki's approach consists of three steps. 
The first step \cite{NPT} is to construct solutions of the
qKZ equation by exploiting hypergeometric integrals (cf. Section 5 for
the explicit formulas). The solutions are $\mathfrak{sl}_2$-singular 
vectors. They are parameterized by polynomials
$P_{n,l}(X_1,\ldots,X_l|z_1,\ldots,z_n)$ in two sets 
of variables, $X_1,\ldots,X_l$ and $z_1,\ldots,z_n$. 
The variables $X_a$ are related
to the integration variables $\alpha_a$ by 
$X_a=e^{-\alpha_a}$, and the 
variables $z_j$ to $\beta_j$ by $z_j=e^{\beta_j}$.
The polynomial $P_{n,l}$ is skew-symmetric in $X_1,\ldots,X_l$ of degree
less than or equal to $n-1$ in each $X_a$,
and is symmetric in $z_1,\ldots,z_n$.
We call $P_{n,l}$ a deformed cycle. Among them there are null cycles,
i.e., those 
which give rise to vanishing integrals. 
The problem of finding all null cycles 
has been solved by Smirnov \cite{S2} and Tarasov \cite{T}.

The second step \cite{N} is to characterize deformed cycles
which give rise to minimal form factors in the following sense. 
A form factor ${\bf f}=(f_n)_{n\ge 0}$ 
is called $N$-minimal if $f_n=0$ for all $n<N$. 
If this is the case, then the function $f_N$ satisfies 
the zero residue condition, 
\begin{equation}\label{NORES}
{\rm res}_{\beta_N=\beta_{N-1}+\pi i}f_N=0.
\end{equation}
We say that a deformed cycle $P_{N,l}$ is $N$-minimal 
if the corresponding function $f_N$ given by the 
hypergeometric integral
satisfies Axiom 1, Axiom 2 and (\ref{NORES}).
Nakayashiki observed that 
a sufficient condition 
for $P_{N,l}$ to be minimal is given by 
\begin{equation}\label{MINI}
P_{N,l}(X_1,\ldots,X_{l-1},z^{-1}|z_1,\ldots,z_{N-2},z,-z)=0.
\end{equation}
He assumed that (modulo null cycles)
the condition (\ref{MINI}) is also necessary, 
and constructed an explicit basis of the 
space of deformed cycles satisfying it.

The third step \cite{N2} is to construct a tower ${\bf f}$,
starting from $f_N$ corresponding to each of 
the basis elements of minimal cycles. 
If $f_n$ and $f_{n+2}$ are given by deformed cycles
$P_{n,l}$ and $P_{n+2,l+1}$, respectively, 
then Axiom 3 (with $n$ replaced by $n+2$) is valid if 
\bea
&&P_{n+2,l+1}(X_1,\ldots,X_l,z^{-1}|z_1,\ldots,z_n,z,-z)
\nn\\
&&\hskip10pt=
z^{-n-1}\prod_{a=1}^l(1-X_a^2z^2)\cdot
P_{n,l}(X_1,\ldots,X_l|z_1,\ldots,z_n).
\label{LINK1}
\ena
In \cite{N2}, Nakayashiki constructed directly 
a tower of deformed cycles 
$\{P_{n,l}\}$
satisfying
the linking condition (\ref{LINK1}),
starting from each minimal cycle he has constructed in \cite{N}.

In \cite{JMT}, a different approach was taken in the construction of
minimal cycles. In that paper, 
the restricted sine-Gordon model (RSG) was considered.
The hypergeometric integrals for ITM and RSG 
have different form, 
but the deformed cycles are exactly the same.
Let us denote the space of minimal deformed 
cycles with fixed $N$ by 
$W_N=\oplus_{0\leq l\leq N}W_{N,l}$. In \cite{JMT}, the space $W_N$ was
constructed as a representation space of a subalgebra of
$U_{\sqrt{-1}}(\widetilde{\mathfrak{sl}}_2)$ 
over the ring $R_N$ of symmetric polynomials in $z_1,\ldots,z_N$, and
it was shown that the total space $W_N$ is created from 
the constant polynomial $1_N\in W_{N,0}$ 
by the $R_N$-linear action of this subalgebra.

In the present paper we continue the study of form factors along
the same direction. We generalize the result in several respects. 
First, we consider deformed cycles $P_{n,l}$ 
which are symmetric Laurent polynomials in the variables $z_1$,$\ldots$,$z_n$.
Namely, we consider both chiralities simultaneously. 
We denote the space of deformed cycles in this extended sense by $\Ch_{n,l}$, 
and that of minimal deformed cycles by $\widehat W_{n,l}$.
Second, we consider the action of the full quantum algebra 
$U_{\sqrt{-1}}(\widetilde{\mathfrak{sl}}_2)$ on 
$\widehat W_n=\oplus_{0\leq l\leq n}\widehat W_{n,l}$.
By doing so, we no longer need to introduce multiplication 
by symmetric Laurent polynomials `by hand',
since it is incorporated as a part of this action on the subspace $\Wh_{n,0}$. 
Finally, and most importantly, we consider also the 
action of $U_{\sqrt{-1}}(\widetilde{\mathfrak{sl}}_2)$
on towers of polynomials. 

Let us elucidate the last point. We are interested in 
a tower of deformed cycles ${\bf p}=(P_{n,l})_{n-2l=m}$ 
satisfying the linking condition 
(\ref{LINK1}), 
where $m$ is a fixed integer. 
As we have seen, such a sequence gives rise to a form factor
satisfying Axiom 3 as well. 
We will refer to ${\bf p}$ as an $\infty$-cycle of weight $m$. 
(
For the precise definition, see Section \ref{subsec:links-of-cycles}.)
The following is an example of $\infty$-cycles of weight $m$, 
\[
{\bf 1}_m=(1_m,X^{m+1},X^{m+1}\wedge X^{m+3},X^{m+1}\wedge
X^{m+3}\wedge X^{m+5},\ldots).
\]
Here we have used the wedge notation for skew-symmetric polynomials
(see \eqref{eq:def-wedge} below).
The action of the quantum algebra extends to 
the space $\widehat C_{n,l}\otimes\C(z_1,\cdots,z_n)$
of polynomials whose coefficients are rational functions 
in $z_1,\ldots,z_n$. 
Therefore it acts naturally on sequences of polynomials 
componentwise. 
Consider the orbit of the particular $\infty$-cycles
given above with $m=0,1$, 
\bea
\Zhc:=\U(\sltt).{\bf 1}_0+\U(\sltt).{\bf 1}_1. 
\label{ztc}
\ena
We show that any element in this space is 
an $\infty$-cycle, that is, 
a sequence of deformed cycles (in particular, they are 
symmetric Laurent polynomials in $z_1,\cdots,z_n$),  
satisfying the linking condition (\ref{LINK1}).
The space \eqref{ztc} is filtered by submodules 
$\Zhc_N$ consisting of $N$-minimal $\infty$-cycles 
${\bf p}=(P_{n,l})$ with $P_{n,l}=0$ for $n<N$. 
We show that the natural map
\be
\oplus_{N\ge 0}\Zhc_N/\Zhc_{N+1}
\longrightarrow \oplus_{N\ge 0}\Wh_N
\en
is an isomorphism of $\U(\sltt)$-modules. 
In particular, any minimal deformed cycle can be lifted
to an $\infty$-cycle, and hence gives rise to a form factor. 
This gives an alternative proof of Nakayashiki's result 
\cite{N2} and extends it in the presence of both chiralities.
These are the main results of the present paper. 


As was shown in \cite{N}, 
the character of the space  
of even, polynomial $\infty$-cycles, i.e.,
$\oplus_{N:{\rm even}}W_N$, 
coincides with the character of the 
level $1$ 
basic representation for $\widehat{\mathfrak{sl}}_2$.
In our convention, it is more natural to think of this 
character as the level $(-1)$-character. 
When we introduce negative powers in $z_j$, 
the character becomes ill defined.   
Instead of dealing with the full space
$\Zhc$, 
we fix a non-negative integer $L\in \Z_{\ge 0}$ 
and consider 
the subspace consisting 
of $N$-minimal deformed cycles $P_{N,l}$ 
such 
that 
$(z_1\cdots z_N)^{L}P_{N,l}$ 
are polynomials in $z_1,\ldots,z_N$. 
Then the sum of the corresponding characters 
(over $N\equiv i\bmod 2$ with $i=0,1$ fixed) has  
the product form $\chi_i(q^{-1},z)\cdot \chi_0(q,z;L)$, 
where $\chi_i(q^{-1},z)$ is the level $(-1)$-character  
and $\chi_0(q,z;L)$ denotes the character of a 
Demazure subspace of the level $1$ irreducible 
module with highest weight $\Lambda_0$. 
The latter tends to $\chi_0(q,z)$ as $L\to\infty$.  
This leads us to conjecture that 
$\Zhc$ is isomorphic 
as a $U_{\sqrt{-1}}(\widetilde{\mathfrak{sl}}_2)$-module 
to the tensor procuct 
of level $(-1)$- and level $1$-modules.
We plan to address this point 
in our next paper. 

Let us give some remarks. 
First, 
the space of $\infty$-cycles is 
not the same as that of form factors. 
On one hand we should 
take into account the quotients by null cycles, 
and on the other hand we should 
incorporate form factors other than the singular vectors
with respect to the ${\mathfrak{sl}}_2$-action.
At the level of characters, these two effects cancel each other. 
The problem of determining 
the symmetries of the space of form factors themselves
is beyond the scope of the present paper. 

Second, lifting of a minimal form factor $f_N$ into a tower  
is not unique.  
Starting from a given lifting, we can add infinitely many $\infty$-cycles
that are minimal and of increasing degrees of minimality,
without changing the degree $N$ part $f_N$. 
Note, however, that such an infinite sum of $\infty$-cycles is not an $\infty$-cycle
in our definition of $\Zhc$. It is not clear if Nakayashiki's extension 
belongs to our space $\Zhc$ (though it is clear that it belongs to
the completion of $\Zhc$). 
Since form factors are determined only up to null cycles, 
the identification of form factors is not a simple problem.
We illustrate this point by some examples. 
Form factors corresponding to the $su(2)$ currents 
given in \cite{S1} belong to the space $\Zhc$.
This is true at the level of $\infty$-cycles. On the other hand,
the known formula \cite{S1} for form factors corresponding to 
the energy-momentum tensor 
arise from $\Zhc$, but only modulo null cycles.



The plan of the paper is as follows. In Section \ref{sec:deformed-cycle}, 
we introduce our notation on the quantum loop algebra at roots of unity, 
and give an action of $\U(\sltt)$ on the space of 
deformed cycles. 
In Section \ref{infinite-cycle} 
we introduce the space of $\infty$-cycles. 
The linking 
condition is preserved by the action of
$U_{\sqrt{-1}}(\widetilde{\mathfrak{sl}}_2)$.
The main result is stated in Theorem 
\ref{thm:main-result}. 
We also calculate the character of the truncated space 
in terms of the level one irreducible characters and the Demazure characters. 
In Section \ref{formfactor}, we apply the results on the $\infty$-cycles
to the form factors. Some technical matters are given in Appendices.
Appendix \ref{app:1} is devoted to the derivation of
the $U_{\sqrt{-1}}(\widetilde{\mathfrak{sl}}_2)$ 
action on the tensor product of evaluation modules. In Appendix \ref{app:2}
we give a proof of a lemma regarding certain null cycles.

\section{The space of deformed cycles}
\label{sec:deformed-cycle}
\subsection{Quantum loop algebra}\label{subsec:loop}

In this subsection, we recall some basic facts
about $U_q(\widetilde{\mathfrak{sl}}_2)$ at roots of unity.
Our basic reference is \cite{CP}. 

Let $\C(q)$ be the field of rational functions 
in indeterminate $q$. 
The quantum loop algebra 
$U_q(\widetilde{\mathfrak{sl}}_2)$ 
is a $\C(q)$-algebra generated by
$x^{\pm}_k$ ($k\in\Z$), $a_n$ ($n\in\Z\backslash\{0\}$) 
and $t_1^{\pm 1}$, with the defining relations
\bea
&&[t_1,a_n]=0,\quad [a_m,a_n]=0,
\label{Dr1}\\
&&
t_1x^{\pm}_kt_1^{-1}=q^{\pm 2}x^{\pm}_k,
\label{Dr2}\\
&&[a_n,x_k^{\pm}]=\pm\frac{[2n]}{n}x^{\pm}_{k+n},
\label{Dr3}\\
&&x^{\pm}_{k+1}x^{\pm}_l-q^{\pm 2}x_l^{\pm}x_{k+1}^\pm
=q^{\pm 2}x^{\pm}_{k}x^{\pm}_{l+1}-x_{l+1}^{\pm}x_{k}^\pm,
\label{Dr4}\\
&&[x^+_k,x_l^-]=
\frac{\varphi^+_{k+l}-\varphi^{-}_{k+l}}{q-q^{-1}}.
\label{Dr5}
\ena
Here 
\be
\sum_{k\in \Z}\varphi^{\pm}_{\pm k}z^k
=t_1^{\pm 1}\exp\left(\pm(q-q^{-1})\sum_{n=1}^\infty
a_{\pm n}z^n\right),
\en
and $[j]=(q^j-q^{-j})/(q-q^{-1})$.
We use the notation
\be
&&x^{(n)}=\frac{x^n}{[n]!},
\quad 
[n]!=\prod_{j=1}^n[j].
\en
Let ${\Ures}^{\pm}$ be the $\C[q,q^{-1}]$-subalgebra 
of $U_q(\widetilde{\mathfrak{sl}}_2)$ 
generated by $(x_n^{\pm})^{(r)}$ ($n\in\Z$, $r\in\Z_{\ge 0}$). 
Let ${\Ures}^0$ be the $\C[q,q^{-1}]$-subalgebra 
generated by 
$t_1^{\pm 1}$, $\qbin{t_1;n}{r}$ ($n\in\Z$, $r\in\Z_{\ge 0}$) 
and 
$\tilde{a}_n$ ($n\in \Z\backslash\{0\}$), where 
\be
\tilde{a}_n=\frac{n}{[n]}a_n,
\quad
\qbin{t_1;n}{r}=
\prod_{s=1}^r
\frac{t_1q^{n-s+1}-t_1^{-1}q^{-n+s-1}}{q^s-q^{-s}}.
\en
Let further $\Ures$ be the $\C[q,q^{-1}]$-subalgebra 
generated by ${\Ures}^\pm$ and ${\Ures}^0$. 
We have the triangular decomposition 
(\cite{CP}, Proposition 6.1)
\be
\Ures={\Ures}^-\cdot{\Ures}^0\cdot{\Ures}^+.
\en
We will use also $\C[q,q^{-1}]$-subalgebras 
$B^\pm_q$ generated by the following elements:
\be
&&B^+_q~:~(x^+_{n})^{(r)},~~
(x^-_{n})^{(r)},~~\tilde{a}_m,~~t_1^{\pm 1},
\\
&&B^-_q~:~
(x^+_{-n})^{(r)},~~
(x^-_{-n})^{(r)},~~\tilde{a}_{-m},~~t_1^{\pm 1},
\en
where $n,r$ run over $\Z_{\ge 0}$
and $m$ over $\Z_{>0}$, respectively. 

Introduce the generating series
\be
&&
\Xc^+_{\ge 0}(t):=\sum_{n\ge 0}x^+_n(q^{-1}t)^n,
\quad
\Xc^+_{< 0}(t):=\sum_{n< 0}x^+_n(q^{-1}t)^n,
\\
&&
\Xc^-_{> 0}(t):=\sum_{n> 0}x^-_n(q^{-1}t)^n,
\quad
\Xc^-_{\le 0}(t):=\sum_{n\le 0}x^-_n(q^{-1}t)^n.
\en
{}By Corollary 4.6 in \cite {CP} and the Remark below it, 
${\Ures}^+$ (resp. ${\Ures}^-$)
is also generated over $\C[q,q^{-1}]$ 
by the coefficients of 
$(\Xc^+_{\ge 0}(t))^{(r)}$ and $(\Xc^+_{< 0}(t))^{(r)}$ 
(resp. $(\Xc^-_{> 0}(t))^{(r)}$ and $(\Xc^-_{\le 0}(t))^{(r)}$), 
where $r$ runs over $\Z_{\ge 0}$.

For any non-zero complex number $\epsilon$, 
the specialization $U_\epsilon$ is defined to be 
$\Ures\otimes_{\C[q,q^{-1}]}\C$ by the ring homomorphism
$\C[q,q^{-1}]\rightarrow\C$ sending $q$ to $\epsilon$.
The subalgebras $U_\epsilon^\pm={\Ures}^\pm$, 
$U^0_\epsilon={\Ures}^0$ and $B^\pm_\epsilon$ 
are defined similarly. 

In this paper we restrict to the case $\U$. 
Define in $\Ures$ 
\be
a_\pm(t):=\sum_{\pm n>0}\tilde a_nt^n,
\en
and denote their images in $\U$ by the same symbols. 
The subalgebras $\U^{\pm}$, $\U^0$ and $B^\pm_{\sqrt{-1}}$ 
are generated over $\C$ by 
(the coefficients of) the following elements:
\bea
\U^-&:&\Xc^-_{>0}(t),\quad \Xc^-_{>0}(t)^{(2)},
\quad
\Xc^-_{\le 0}(t),\quad \Xc^-_{\le 0}(t)^{(2)},
\label{Um}\\
\U^0&:&
t_1^{\pm 1},\quad \qbin{t_1;n}{2}~~(n\in\Z),
\quad 
a_{\pm}(t),
\label{U0}\\
\U^+&:&\Xc^+_{\ge0}(t),\quad \Xc^+_{\ge0}(t)^{(2)},
\quad
\Xc^+_{<0}(t),\quad \Xc^+_{<0}(t)^{(2)},
\label{Up}
\\
B^+_{\sqrt{-1}}&:&
\Xc^+_{\ge0}(t),\quad \Xc^+_{\ge0}(t)^{(2)},
\quad \Xc^-_{>0}(t),\quad \Xc^-_{>0}(t)^{(2)},
\label{Bp}
\\
&&\quad x_0^-,\quad (x_0^-)^{(2)},
\quad a_+(t),\quad t_1^{\pm 1},
\nn
\\
B^-_{\sqrt{-1}}&:&
\Xc^+_{<0}(t),\quad \Xc^+_{<0}(t)^{(2)},
\quad \Xc^-_{\le 0}(t),\quad \Xc^-_{\le 0}(t)^{(2)},
\label{Bm}
\\
&&\quad x_0^+,\quad (x_0^+)^{(2)},
\quad a_-(t),\quad t_1^{\pm 1}.
\nn
\ena

We assign the degree and weight to $\U$ as follows.
\be
&&\deg (x^\pm_n)^{(r)}=nr,~~\deg \tilde{a}_m=m,
~~\deg t_1^{\pm 1}=0,\\
&&
\wt (x^\pm_n)^{(r)}=\pm 2r,~~
\wt \tilde{a}_m=0,~~
~~\wt t_1^{\pm 1}=0. 
\en

\subsection{Representation of $\U$ on the space of
polynomials}\label{action-on-cycle}

Until the end of this section, 
we fix a non-negative integer $n$.
Let $K_n=\C(z_1,\cdots,z_n)$ be the field of rational functions
in $z_{1}, \ldots , z_{n}$, and let $\widehat R_n$ be its subring
consisting of symmetric Laurent polynomials.
Consider the vector space over $K_n$
\be
&&
A_n:=\oplus_{l=0}^nA_{n,l},\quad 
A_{n,l}:=\wedge^l\bigl(\oplus_{j=0}^{n-1} K_nX^j\bigr).
\en
By definition we set $A_{n, l}=0$ if $l<0$. 
Note that $A_{n, 0}=K_{n}$ and 
$A_{n, l}=0$ if $l>n$. 
For $0\le l\le n$, we identify an element $P\in A_{n,l}$ with 
a skew-symmetric polynomial in the variables
$X_1,\ldots,X_l$ with coefficients in $K_n$, 
of degree at most $n-1$ in each $X_j$. 
We use the wedge product notation 
for $P_1(X_1,\cdots,X_{l_1})\in A_{n,l_1}$
and $P_2(X_1,\cdots,X_{l_2})\in A_{n,l_2}$, 
\bea
&& 
P_1\wedge P_2(X_1,\cdots,X_{l_1+l_2}) 
\label{eq:def-wedge} \\ 
&& 
:=\frac{1}{l_1!l_2!}{\rm Skew}\,
\left( 
P_1(X_1,\cdots,X_{l_1})P_2(X_{l_1+1},\cdots,X_{l_1+l_2}) 
\right), \nn
\ena
where ${\rm Skew}$ stands for the skew-symmetrization 
\be
{\rm Skew}\,f(X_1,\cdots,X_m):=
\sum_{\sigma\in \mathfrak{S}_m}({\rm sgn}\,\sigma)
f(X_{\sigma(1)},\cdots,X_{\sigma(m)}). 
\en

We shall introduce a $K_n$-linear action of $\U$ on $A_n$. 
For that purpose, let us prepare some notation.
Set 
\be
&&\Theta_n(X):=\prod_{j=1}^n(1-z_jX),
\\
&&\Theta_n(X_1,X_2):=
\Theta_n(X_1)\Theta_n(X_2)-\Theta_n(-X_1)\Theta_n(-X_2).
\en
Define also 
\be
&&F_n(t|X):=\frac{t}{2(X-t)}\frac{\Theta_n(t,-X)}{\Theta_n(t)},
\\
&&F^{(2)}_{n}(t|X_1,X_2)
:=\frac{t}{X_1+t}\frac{t}{X_2+t}
\frac{X_1-X_2}{X_1+X_2}\Theta_n(X_1,X_2)\\
&&
+\frac{t}{X_1-t}\frac{t}{X_2+t}
\frac{\Theta_n(-X_2)}{\Theta_n(t)}\Theta_n(t,-X_1)
{}-\frac{t}{X_2-t}\frac{t}{X_1+t}
\frac{\Theta_n(-X_1)}{\Theta_n(t)}\Theta_n(t,-X_2).
\en
$F_n(t|X)$ is a polynomial in $X$ of degree $n-1$ because $\Theta_n(t,-X)$
is divisible by $X-t$.
At $t^{\pm 1}=0$, 
it has a power series expansion in $t^{\pm 1}$
whose coefficients 
are symmetric polynomials in $z_1^{\pm 1},\cdots,z_n^{\pm 1}$. 
$F^{(2)}_{n}(t|X_1,X_2)$ has similar properties.

Let us define the action of the generators 
$g\in \U$ on $P\in A_{n,l}$. 
For $l=0$ or $1$ and $g\in \U^+$, we set
\be
&&\Xc^{+}_{\ge 0}(t).P=0,\quad
\Xc^{+}_{< 0}(t).P=0\quad (l=0),
\\
&&\Xc^{+}_{\ge 0}(t)^{(2)}.P=0,\quad
\Xc^{+}_{< 0}(t)^{(2)}.P=0\quad (l=0,1).
\en
In the other cases, we define the action as follows.
\bea
&&t_1^{\pm 1}.P:=i^{\pm(n-2l)}.P,
\label{act0}\\
&&\Xc^-_{>0}(t). P:=F_n(t)\wedge P,
\label{act1}\\
&&-4i\Xc^-_{>0}(t)^{(2)}. P:=F^{(2)}_n(t)\wedge P,
\label{act2}\\
&&a_+(t).P:=\sum_{j=1}^n\frac{z_jt}{1-z_jt}P(X_1,\cdots,X_l)\label{act4}\\
&&+\sum_{p=1}^l\Bigl\{\frac{t}{X_p-t}\Bigl(\frac{\Theta_n(X_p)}{\Theta_n(t)}
P(X_1,\cdots,\overset{\mbox{\tiny $\overset{p}{\smallsmile}$}}{t},\cdots,X_l)
{}-P(X_1,\ldots,X_l)\Bigr)\nn\\
&&\hskip100pt+(t\rightarrow-t)\Bigr\},\nn\\
&&a_-(t).P:=-\sum_{j=1}^n\frac1{1-z_jt}P(X_1,\cdots,X_l)\label{act5}\\
&&-\sum_{p=1}^l\Bigl\{\frac{t}{X_p-t}\Bigl(\frac{\Theta_n(X_p)}{\Theta_n(t)}
P(X_1,\cdots,\overset{\mbox{\tiny $\overset{p}{\smallsmile}$}}{t},\cdots,X_l)
{}-\frac {X_p}tP(X_1,\ldots,X_l)\Bigr)\nn\\
&&\hskip100pt+(t\rightarrow-t)\Bigr\},\nn\\
&&i^{1-n}\Xc^+_{\ge 0}(t). P
:=\frac{1}{\Theta_n(t)}P(X_1,\cdots,X_{l-1},t),
\label{act6}\\
&&i(-1)^{n+1}\Xc^+_{\ge 0}(t)^{(2)}. P
:=\sum_{a=1}^n\res_{u=z_a^{-1}}
\left(\frac{P(X_1,\cdots,X_{l-2},-u,u)}
{\Theta_n(-u)\Theta_n(u)}\frac{du}{u-t}
\right).
\label{act7}
\ena
In \eqref{act1},\eqref{act2},\eqref{act4},\eqref{act6}
and \eqref{act7},
the right hand sides stand for the power series expansion 
in $t$, and in \eqref{act5} the expansion in $t^{-1}$. 
Note that the factor $X_p\pm t$ in the denominator of (\ref{act4})
or (\ref{act5}) divides the numerator.
Define also $-\Xc^-_{\le 0}(t).P$, 
$-4i\Xc^-_{\le 0}(t)^{(2)}.P$, 
$-i^{1-n}\Xc^+_{<0}(t).P$ and 
$i(-1)^{n+1}\Xc^+_{<0}(t)^{(2)}.P$ 
by the expansion at $t=\infty$ of the right hand side of 
\eqref{act1}, \eqref{act2}, \eqref{act6}, \eqref{act7}, 
respectively. 
Note, in particular, that
\begin{equation}
\tilde a_m.P:=\bigl(\sum_{j=1}^nz_j^m\bigr)P
\qquad (\mbox{$m\in\Z$ : odd}). 
\label{act3}
\end{equation}

\begin{prop}\label{prop:Umod}
With the above rule, $A_n$ is a $\U$-module. 
\end{prop}
This action of $\U$ is essentially the one on 
the $n$-fold tensor product of 
two-dimensional evaluation modules in disguise. 
We give a proof of Proposition \ref{prop:Umod} in Appendix A. 

The following lemma will be used later. 
\begin{lem}\label{lem:action-of-an}
On the subspace $A_{n,0}=K_n$, we have
\be
\tilde{a}_m.P=\bigl(\sum_{j=1}^nz_j^m\bigr).P
\qquad (m\in\Z\backslash\{0\},P\in A_{n,0}).
\en
We have $\Rh_n=\U^0.1_n$, where $1_n\in\Rh_n$ denotes the unit.
\end{lem}
\begin{proof}
This is an immediate consequence of 
\eqref{act4}--\eqref{act5} and \eqref{act3}.
\end{proof}

\subsection{Deformed cycles}\label{subsec:deformedcycle}

For the application to form factors, 
we introduce several submodules of $A_n$ 
\be
A_n\supset \Ch_n\supset \Dh_n\supset \Wh_n
\en
to be defined as follows. 

Set
\be
&&
\Ch_n:=\oplus_{l=0}^n\Ch_{n,l},\quad 
\Ch_{n,l}:=\wedge^l\bigl(\oplus_{j=0}^{n-1} 
\hat{R}_nX^j\bigr).
\en
An element of $\Ch_n$ will be referred to as a deformed cycle.
We say that $P\in\Ch_{n,l}$ is weakly minimal if
\bea
&&
\mbox{$P=0$ if $X_{l-1}^{-1}=-X_l^{-1}=z_{n-1}=-z_n$},
\label{w-minimal}
\ena
and minimal if 
\bea
&&
\mbox{$P=0$ if $X_l^{-1}=z_{n-1}=-z_n$}.
\label{minimal}
\ena
These conditions arise in the study of form factors
(see Lemma \ref{TRIVLEM} below). 
Denote by $\Dh_{n,l}$ (resp. $\Wh_{n,l}$) 
the subspace of weakly minimal (resp. minimal) 
elements of $\Ch_{n,l}$. 
We set 
\be
\Dh_n:=\oplus_{l=0}^n\Dh_{n,l},~~
\Wh_n:=\oplus_{l=0}^n\Wh_{n,l}.
\en
Let further $R_n$ be the subring of $\Rh_n$ 
consisting of symmetric polynomials.
Replacing $\widehat{R}_n$ by $R_n$
in the above, we define the spaces 
$C_n=\oplus_{l=0}^nC_{n,l}$, 
$D_n=\oplus_{l=0}^nD_{n,l}$, 
$W_n=\oplus_{l=0}^nW_{n,l}$.

\begin{lem}\label{lem:nopole1}
\begin{enumerate}
\item The spaces $\Dh_n$, $\Wh_n$ are $\U$-submodules. 
\item 
The spaces $D_n$, $W_n$ are $B^+_{\sqrt{-1}}$-submodules. 
\end{enumerate}
\end{lem}
\begin{proof}
{}From the formulas \eqref{act0}--\eqref{act7}, 
it is clear that $\Dh_n$, $\Wh_n$ are stable under the action of
all generators listed in \eqref{Um}, \eqref{U0} and \eqref{Up}, except for
$\Xc^+_{\ge 0}(t)^{(2)}$, $\Xc^+_{< 0}(t)^{(2)}$. 
The only subtle point about the latter elements 
is that, when acted on $P\in\Ch_n$, they give rise to 
symmetric rational functions in $z_1,\cdots,z_n$ 
which have a pole on $z_a+z_b=0$ in general.
We show that for elements $P\in \Dh_n$ 
these poles do not appear, 
and that the (weak-) minimality condition is preserved. 

As an example, 
consider the case $Q=i(-1)^n\Xc^+_{\ge 0}(t)^{(2)}.P$, 
$P\in \Dh_{n,l}$ with $l\ge 2$.  
Explicitly we have 
\be
&&Q(X_1,\cdots,X_{l-2}|z_1,\cdots,z_n)
=\frac{1}{2}\sum_{a=1}^n\frac{1}{\prod_{b(\neq a)}(1-z_b^2/z_a^2)}
\frac{1}{1-z_at}
\\
&&\qquad \quad
\times P(X_1,\cdots,X_{l-2},z_a^{-1}, -z_a^{-1}|z_1,\cdots,z_n).
\en
Since the right hand side is symmetric in $z_1,\cdots,z_n$, 
the only possible poles are those at $z_a+z_b=0$ ($a\neq b$). 
However, those poles are absent because of \eqref{w-minimal}.
Let us verify that $Q\in\Dh_{n,l-2}[[t]]$. 
If $l<4$, there is nothing to show. 
Suppose $l\ge 4$. 
Setting $z_{n-1}=z$ and $z_n=-z$ we find 
\be
&&Q(X_1,\cdots,X_{l-2}|z_1,\cdots,z_{n-2},z,-z)
\\
&&=\frac{1}{2}
\sum_{a=1}^{n-2}
\frac{P(X_1,\cdots,X_{l-2},z_a^{-1},-z_a^{-1}|z_1,\cdots,z_{n-2},z,-z)}
{(1-tz_a)(1-z^2/z_a^2)^2\prod_{1\le b\le n-2\atop b(\neq a)}(1-z_b^2/z_a^2)}
\\
&&+\frac{1}{4}
\frac{1}{\prod_{b=1}^{n-2}(1-z_b^2/z^2)}
\Bigl(
\frac{z}{1+tz}\frac{\partial P}{\partial z_{n-1}}
(X_1,\cdots,X_{l-2},z^{-1},-z^{-1}|z_1,\cdots,z_{n-2},z,-z)
\\
&&\qquad+\frac{z}{1-tz}\frac{\partial P}{\partial z_{n}}
(X_1,\cdots,X_{l-2},z^{-1},-z^{-1}|z_1,\cdots,z_{n-2},z,-z)
\Bigr).
\en
Under further specialization $X_{l-3}=-X_{l-2}=z^{-1}$, 
the first term vanishes by \eqref{w-minimal}, 
and the rest is $0$ by skew symmetry.

In the same way, \eqref{minimal} is also preserved. 
The case $\Xc^+_{<0}(t)^{(2)}$ 
and the remaining assertions can be shown similarly.
\end{proof}

\begin{prop}\label{prop:W-is-cyclic}
Let $1_n\in \Ch_n$ be the unit element. Then we have
\be
&&\Wh_n=\U.1_{n},
\quad W_n=B^+_{\sqrt{-1}}.1_{n}. 
\en
\end{prop}
\begin{proof}
By Lemma \ref{lem:nopole1}, $\U$ acts on $\Wh_n$. 
It is known \cite{JMT} that $W_n$ is generated from 
$1_n$ by the action of $x_0^-$, $(x_0^-)^{(2)}$, 
$\Xc^-_{>0}(t)$, $\Xc^-_{>0}(t)^{(2)}$ and the 
multiplication by elements of $R_n$. 
By Lemma \ref{lem:action-of-an}, 
$R_n=(B^+_{\sqrt{-1}}\cap \U^0).1_n$.
Since $\Xc^+_{\ge0}(t)$, $\Xc^+_{\ge0}(t)^{(2)}$ kill $1_n$, 
and since the action of $\U$ is $R_n$-linear, 
we obtain $W_n=B^+_{\sqrt{-1}}.1_{n}$. 
The other assertion follows in a similar manner.
\end{proof}

Consider the tensor product of evaluation representation 
$\rho_{z_1}\otimes\cdots\otimes\rho_{z_n}$ specialized to 
$q=\sqrt{-1}$ (see \eqref{eq:eval}).
As explained in Appendix \ref{app:1}, 
the action of $\U$ on $\Wh_n$ is induced from 
$\rho_{z_1}\otimes\cdots\otimes\rho_{z_n}$ under
the identification of
$1_n$ with $v_+^{\otimes n}$. 
Proposition \ref{prop:W-is-cyclic} implies that 
\begin{cor}
As $\U$-module, $\Wh_n$ is isomorphic to the subrepresentation of 
$\rho_{z_1}\otimes\cdots\otimes\rho_{z_n}$ 
generated by $v_+^{\otimes n}$. 
\end{cor}

Define the degree on $\widehat{C}_{n, l}$ by 
\begin{equation}
\deg_{0}X_{a}=-1, \quad \deg_{0}z_{j}=1. \label{DEG0}
\end{equation}
We also assign a weight $m=n-2l$ to $P\in \Ch_{n,l}$.
With this definition, 
$\Dh_n$ and $\Wh_n$ are bi-graded $\U$-modules. 


\section{Action of $\U(\sltt)$ 
on the space of $\infty$-cycles}
\label{infinite-cycle}

In this section we introduce certain sequences of cycles 
which we call $\infty$-cycles, and 
define an action of $\U$ on them.
\subsection{Links of cycles}\label{subsec:links-of-cycles}

Let $P_{n,l}\in \Ch_{n,l}$, $P_{n+2,l+1}\in \Ch_{n+2,l+1}$. 
We say that the pair $(P_{n,l},P_{n+2,l+1})$ is a {\it link} if 
\bea\label{link-def}
&&P_{n+2,l+1}(X_1,\cdots,X_l,z^{-1}|z_1,\cdots,z_n,z,-z)
\nn\\
&&=z^{-n-1}\prod_{a=1}^l(1-X_a^2z^2)\cdot
P_{n,l}(X_1,\cdots,X_l|z_1,\cdots,z_n).
\ena
The condition \eqref{link-def} appeared in \cite{NT} in the following 
equivalent form. 
\begin{lem}\label{P*}
A pair $(P_{n,l},P_{n+2,l+1})$ is a link if and only
if there exists a 
\be
\Ps_{n,l}(X_1,\cdots,X_l|X_{l+1}|z_1,\cdots,z_{n}|z)
\en 
with the following properties.
\begin{enumerate}
\item $\Ps_{n,l}$ 
is a skew-symmetric polynomial in $X_1,\cdots,X_{l}$
with $\deg_{X_i}\Ps_{n,l}\le n-1$ ($1\le i\le l$). 
It is a polynomial in $X_{l+1}$ with 
$\deg_{X_{l+1}}\Ps_{n,l}\le n+1$. 
\item $\Ps_{n,l}$ is a symmetric Laurent polynomial in $z_1,\cdots,z_n$, and 
an even Laurent polynomial in $z$.
\item We have
\bea
&&P_{n+2,l+1}(X_1,\cdots,X_{l+1}|z_1,\cdots,z_n,z,-z)
\label{recursion-link}\\
&&=\frac{1}{l!}{\rm Skew}\left( 
\prod_{1\le a\le l}(1-X_a^2z^2)\,
\Ps_{n,l}(X_1,\cdots,X_l|X_{l+1}|z_1,\cdots,z_n|z) \right),
\nn\\
&&\Ps_{n,l}(X_1,\cdots,X_l|z^{-1}|z_1,\cdots,z_n|z)
=z^{-n-1}P_{n,l}(X_1,\cdots,X_l|z_1,\cdots,z_n).
\ena
Here ${\rm Skew}$ 
in \eqref{recursion-link} stands for the skew-symmetrization 
with respect to $X_{1}, \ldots , X_{l+1}$. 
\end{enumerate}
\end{lem}
\begin{proof}
The `if' part is evident. To see the `only if' part, suppose \eqref{link-def} 
holds. Since 
\be
&&P_{n+2,l+1}(X_1,\cdots,X_{l+1}|z_1,\cdots,z_n,z,-z)\\
&&
-\frac{1}{l!}{\rm Skew}\Bigl(X_{l+1}^{n+1}\prod_{a=1}^n(1-X_a^2z^2)
\cdot P_{n,l}(X_1,\cdots,X_l|z_1,\cdots,z_n)\Bigr)
\en
has a zero at $X_{l+1}=\pm z^{-1}$, it can be written as 
$\prod_{a=1}^{l+1}(1-X_a^2z^2)\cdot Q(X_1,\cdots,X_{l+1}|z_1,\cdots,z_n|z)$. 
Setting 
\be
&&\Ps_{n,l}(X_1,\cdots,X_l|X_{l+1}|z_1,\cdots,z_n|z)
\\
&&=X_{l+1}^{n+1}P_{n,l}(X_1,\cdots,X_l|z_1,\cdots,z_n)
+\frac{1}{l+1}(1-X_{l+1}^2z^2)Q(X_1,\cdots,X_{l+1}|z_1,\cdots,z_n|z),
\en
one easily verifies the conditions mentioned above. 
\end{proof}

The following lemma is obvious from the definition.
\begin{lem}\label{TRIVLEM}
\begin{enumerate}
\item If $(P_{n,l},P_{n+2,l+1})$ is a link, then 
$P_{n,l}$ and $P_{n+2,l+1}$ are weakly minimal. 
\item For $P_{n,l}\in \Ch_n$, 
$(0, P_{n,l})$ is a link if and only if $P_{n,l}$ is minimal.
\end{enumerate}
\end{lem} 

Let us study the action of $\U$ on links. 
\begin{prop}\label{prop:link-link}
Let $P_{n,l}\in\Ch_{n,l}$, $P_{n+2,l+1}\in\Ch_{n+2,l+1}$ 
and $g\in \U$. 
If the pair $(P_{n,l},P_{n+2,l+1})$ is a link, 
then $(g.P_{n,l},g.P_{n+2,l+1})$ is also a link. 
\end{prop}
Using the formulas for the action \eqref{act1}--\eqref{act7},
it is straightforward to verify
the linking condition \eqref{link-def}.
We omit the proof.
\subsection{$\infty$-cycles} 

Let $m$ be an integer. 
We call a sequence of cycles 
\be 
{\bf p}=(P_{n, l})_{n, l \ge 0,n-2l=m}, \quad 
P_{n, l} \in \widehat{C}_{n, l} 
\en 
an {\it $\infty$-cycle} if 
$(P_{n, l}, P_{n+2, l+1})$ is a link for any $n \ge 0$. 
The integer $m$ is called the weight of ${\bf p}$. 
We denote by $\Zbc[m]$ the space of $\infty$-cycles of 
weight $m$, and set $\Zbc=\oplus_{m\in\Z}\Zbc[m]$.  
Often we write the component $P_{n,l}$ as $({\bf p})_{n,l}$.
For an $\infty$-cycle ${\bf p}\in\Zbc[m]$, 
define the action of $g\in \U$ by
\be
g.{\bf p}=\bigl(g.P_{n,l})_{n,l\ge 0,n-2l=m}.
\en
\medskip

\noindent 
{\it Example.}\quad
{}For each $m\in\Z_{\ge 0}$, there is a distinguished 
$\infty$-cycle ${\bf 1}_m\in \Zbc[m]$ given by 
\be 
({\bf 1}_{m})_{m,0}=1_m,
\quad
{\bf 1}_{m+2l,l}=X^{m+1}\wedge\cdots\wedge X^{m+2l-1}
\quad (l\ge 1). 
\en 
{}From the formula \eqref{act1} we find
\be
(i^{-1-m}\Xc^+_{\ge 0}(t).{\bf 1}_m)_{m+2l,l-1}
=t^{m+1}X^{m+3}\wedge\cdots\wedge X^{m+2l-1}+O(t^{m+2}),
\en
and hence 
\be
(-1)^{m+1}x^+_{m+1}{\bf 1}_m={\bf 1}_{m+2},\quad
x^+_n{\bf 1}_m=0\quad (n\le m).\\
\en
In particular, for all $m\ge 0$ we have
\bea
{\bf 1}_m=
\begin{cases}
(-1)^{k}x^+_{2k-1}\cdots x^+_3x^+_1 {\bf 1}_0 & (m=2k),\\
x^+_{2k-2}\cdots x^+_2x^+_0 {\bf 1}_1 & (m=2k-1).\\
\end{cases}
\label{extremal-vectors}
\ena

Next we give a formula for 
$x_{2k-1}^{-} \cdots x_{3}^{-} x_{1}^{-} {\bf 1}_{0} \in \Zbc[2k]$ 
as a nontrivial example. 
For $k=1, 2$ the $\infty$-cycles are given by 
\be 
&& 
(x_{1}^{-}{\bf 1}_{0})_{2+2l, 2+l}= 
i^{1/2}\left( 
\sum_{0 \le a <l+1}e_{2a+1}X^{2a} \right) 
\wedge X \wedge X^{3} \wedge \cdots \wedge X^{2l+1}, \\   
&& 
(x_{3}^{-}x_{1}^{-}{\bf 1}_{0})_{4+2l, 4+l} \\ 
&& {}=  
\left( 
\sum_{0 \le a_{1}<a_{2} < l+2} 
\left| 
\begin{array}{ccc} 
e_{2a_{2}+1} & e_{2a_{2}+2} & e_{2a_{2}+3} \\ 
e_{2a_{1}+1} & e_{2a_{1}+2} & e_{2a_{1}+3} \\ 
0 & 1 & e_{1} 
\end{array} 
\right|  
X^{2a_{1}} \wedge X^{2a_{2}} 
\right) \wedge 
X \wedge X^{3} \wedge \cdots \wedge X^{2l+3},  
\en 
where $e_{a}$ is the $a$-th elementary symmetric polynomial in $z_{j}$'s. 
To write down the formula for general $k$, we set 
\be 
w_{2n, r}^{(0)}&:=&(-1)^{\frac{r(r-1)}{2}} \\ 
&\times& 
\sum_{0 \le a_{1}< \cdots <a_{r} <n}
S_{(2(r-1), \ldots , 2, 0 | 2a_{r}, \ldots , 2a_{1})}(z_{1}, \ldots , z_{2n}) 
X^{2a_{1}} \wedge \cdots \wedge X^{2a_{r}}. 
\en 
Here $S_{(2(r-1), \ldots , 2, 0| 2a_{r}, \ldots , 2a_{1})}(z_{1}, \ldots , z_{2n})$ 
is the Schur polynomial 
associated with the Young diagram given in the Frobenius notation 
by $(2(r-1), \ldots , 2, 0|2a_{r}, \ldots , 2a_{1})$. 
Then the formula is given by 
\be 
(x_{2k-1}^{-} \cdots x_{3}^{-} x_{1}^{-} {\bf 1}_{0})_{2k+2l, 2k+l}= 
i^{k^{2}/2}w_{2k+2l, k}^{(0)} \wedge X \wedge X^{3} \wedge \cdots \wedge X^{2(k+l)-1}. 
\en 
In particular we have 
\be 
(x_{2k-1}^{-} \cdots x_{3}^{-} x_{1}^{-} {\bf 1}_{0})_{2k, 2k}= 
i^{k^{2}/2}\prod_{1 \le j<j' \le 2k}(z_{j}+z_{j'}) \,
X^{0} \wedge X^{1} \wedge  \cdots \wedge X^{2k-1}. 
\en 
It is clear that this cycle is minimal. 
\qed
\medskip

The central object of our study are the following modules.
\be
&&\Zhc:=\Zhc^{(0)}\oplus\Zhc^{(1)}, \quad \Zhc^{(i)}:=\U.{\bf 1}_i
\quad (i=0,1),\\
&&\Zc:=\Zc^{(0)}\oplus\Zc^{(1)}, \quad \Zc^{(i)}:=B^+_{\sqrt{-1}}.{\bf 1}_i
\quad (i=0,1).
\en
We have $\Zhc=\oplus_{m\in\Z}\Zhc[m]$, 
$\Zhc[m]:=\Zhc\cap\Zbc[m]$ and likewise for $\Zc$. 
We also define the degree of an $\infty$-cycle 
${\bf p}=(P_{n, l})$ by 
\bea 
\deg{\bf p}=\frac{n^{2}}{4}+\deg_{0}P_{n, l},
\label{eq:def-degree} 
\ena 
where $\deg_0$ is defined in (\ref{DEG0}).
Since any pair $(P_{n, l}, P_{n+2, l+1})$ in ${\bf p}$ 
is a link, 
the right hand side above is independent of $n$. 

We say an $\infty$-cycle ${\bf p}=(P_{n, l})$ is $N$-{\it minimal} if 
$P_{n, l}=0$ for $n < N$. 
The space $\Zhc$ is filtered by $\U$-submodules
$\Zhc_N$ consisting of $N$-minimal $\infty$-cycles, 
\[
\Zhc=\Zhc_0\supset \Zhc_1\supset\cdots.
\]
Similarly we have a filtration by 
$B^+_{\sqrt{-1}}$-submodules
$\Zc_N=\Zc\cap\Zhc_N$ of $\Zc$. 
If ${\bf p}=(P_{n, l})_{n,l\ge 0,n-2l=m} \in \Zhc_N$, 
then $P_{N, l'}$ is $N$-minimal where $l'=\frac{N-m}{2}$. 
We have therefore natural injective homomorphisms
\bea
\hat{\varphi}&:&
\gr \Zhc:=\oplus_{N\geq0}\Zhc_N/\Zhc_{N+1}\rightarrow 
\oplus_{N\ge 0}\widehat{W}_{N},
\label{isoh}\\
\varphi&:&
\gr \Zc:=\oplus_{N\geq0}\Zc_N/\Zc_{N+1}\rightarrow 
\oplus_{N\ge 0}{W}_{N}, 
\label{iso}
\ena
which respect the bi-grading $(d,m)$ given by $-{\rm deg}$ and weight.
(See Section \ref{section:character} for the reason of the minus sign.)
We are now in a position to state the main result of this paper. 

\begin{thm}\label{thm:main-result}
The map \eqref{isoh} $($resp. \eqref{iso}$)$ is an 
isomorphism of bi-graded $\U$-modules
$($resp. $B^+_{\sqrt{-1}}$-modules$)$.
\end{thm}
\begin{proof}
It suffices to show the surjectivity. 
Take any $P\in\Wh_N$. 
By Proposition \ref{prop:W-is-cyclic},
there exists $g\in\U$ such that $P=g.1_N$. 
{}From \eqref{extremal-vectors}, we can choose 
$g'\in\U$ and $i\in\{0,1\}$ such that ${\bf 1}_N=g'.{\bf 1}_i$. 
Then ${\bf p}=gg'.{\bf 1}_i$ is an $\infty$-cycle which is 
sent to $P$ under the map \eqref{isoh}. 
The case \eqref{iso} is completely parallel. 
\end{proof}
Theorem \ref{thm:main-result}
gives an alternative proof 
of Nakayashiki's result \cite{N2} and extends it 
to the case including negative powers of $z_1,\cdots,z_n$. 
\medskip

\noindent
{\it Remark.}\quad
Fix $m$. For any sequence ${\bf p}_j\in\Zbc[m]$ 
such that $({\bf p}_j)_{n,(n-m)/2}=0$ for $n<j$, 
the infinite sum $\sum_{j\ge 0}{\bf p}_j$ 
is well defined and belongs to $\Zbc[m]$. 
Theorem \ref{thm:main-result} implies
that conversely any 
${\bf p}\in\Zbc[m]$ can be written as such an infinite sum 
with ${\bf p}_j\in\Zhc[m]$. 
\qed

\subsection{Characters}\label{section:character}

In this subsection we study the characters of $\Zc$, $\Zhc$. 
By a character of a bi-graded vector space $V=\oplus_{d,m}V_{d,m}$, 
we mean the generating series
\bea
{\rm ch}_{q,z} V=\sum_{d,m}q^dz^m \dim V_{d,m}. 
\label{def:characterofV}
\ena
In the below we use the bi-grading $(d,m)$ by $-{\rm deg}$ and weight.
In order to match our characters to the irreducible characters
of $\widehat{\mathfrak{sl}}_2$ we put the minus sign for the degree.

\begin{thm}\label{thm:chiral-character}
For $i=0$ or $1$, we have
\begin{equation}
\ch_{q,z}\Zc^{(i)}=\chi_i(q^{-1},z),
\label{CHARAC}
\end{equation}
where
\be
\chi_i(q,z):=\frac{1}{(q)_\infty}\sum_{m\in 2\Z+i}q^{m^2/4}z^m
\en
is the character of the level $1$ integrable 
module with the highest weight $\Lambda_i$ of $\widehat{\mathfrak{sl}}_2$.
\end{thm}
\begin{proof}
This is a consequence of the isomorphism \eqref{iso} and 
Nakayashiki's result \cite{N} on the character of $W_n$. 
\end{proof}
Note that (\ref{CHARAC}) implies that the character of $\Zc^{(i)}$
is equal to that of the level $-1$ integrable module
with the lowest weight $-\Lambda_i$.

On the other hand, the character of $\Zhc^{(i)}$ is ill-defined 
since each weight subspace is infinite dimensional.
To get around this inconvenience, we follow the idea of \cite{Kou} 
and introduce the truncated character. 
For each $L\in\frac12\Z_{\geq0}$, consider the space
\be
\Wh_{N,L}:=\bigl(\prod_{j=1}^Nz_j\bigr)^{-L}\cdot W_N.
\en
We have $\Wh_{N,0}\subset\Wh_{N,1}\subset\cdots$ and
$\Wh_N=\cup_{L\in\Z_{\ge 0}}\Wh_{N,L}$.
We will use below the standard $q$-binomial symbol
\be
\qbin{m}{n}=
\begin{cases}
\displaystyle{\frac{(q)_m}{(q)_n(q)_{m-n}}}
& (m\ge n\ge 0),\\
0 & (\mbox{otherwise}),\\
\end{cases}
\en
where $(z)_n=\prod_{j=0}^{n-1}(1-q^jz)$.
Recall that each integrable $\widehat{\mathfrak{sl}}_2$-module 
has a family of Demazure subspaces
parametrized by the elements of the affine Weyl group. Their characters 
\be
\chi_i(q,z;L)=\sum_{m\in 2\Z+i}
\qbin{\scriptstyle2L}{{\scriptstyle L+\frac m2}}q^{m^2/4}z^m
\quad (i\equiv 2L\bmod 2) 
\en
give polynomial finitizations of the full characters in the sense that 
$\lim_{L\to\infty}\chi_i(q,z;L)=\chi_i(q,z)$. 

\begin{prop}\label{prop:Demazure}
Let $L\in\frac12\Z_{\ge 0}$, $i,j=0,1$ with $j\equiv2L\bmod 2$. 
Then we have 
\bea
{\rm ch}_{q,z}\bigl(\oplus_{N\ge 0\atop N\equiv i \bmod2}
\Wh_{N,L}\bigr)
=
\chi_i(q^{-1},z)\chi_j(q,z;L).
\label{twochirality}
\ena
\end{prop}
\begin{proof}
Without loss of generality, 
we prove \eqref{twochirality} in the summed form over $i=0,1$.
Theorem \ref{thm:chiral-character} is equivalent to the
statement
\be
{\rm ch}_{q^{-1},z}W_{N}=q^{N^2/4}\sum_{l=0}^N\frac{1}{(q)_N}
\qbin{N}{l}z^{N-2l}.
\en
Hence, taking the sum over $i=0,1$, 
the left hand side of \eqref{twochirality} with $q$ replaced by $q^{-1}$
becomes
\bea
\sum_{N\ge l\ge 0}\frac{q^{N^2/4-NL}}{(q)_N}\qbin{N}{l}z^{N-2l}
=\sum_{l,m\in\Z\atop l\ge 0, -m}
z^mq^{m^2/4-Lm}
\frac{q^{lm+l(l-2L)}}{(q)_l(q)_{l+m}}.
\label{trun-char}
\ena
The sum over $l$ can be simplified by 
using the formula (\cite{JMT}, Lemma 6.2)
\bea
\sum_{l\ge 0}\frac{(q^{l+1}z)_\infty}{(q)_l}
q^{l(l-2L)}z^l
=\sum_{s\ge0}\qbin{2L}{s}_{q^{-1}}z^s,
\label{identity}
\ena
wherein ${\scriptstyle \qbin{2L}{s}_{q^{-1}}}$ 
signifies the $q$-binomial symbol
with $q$ replaced by $q^{-1}$. 
Specializing $z=q^m$ in \eqref{identity}, 
we find that the right hand side of 
\eqref{trun-char} becomes
\be
\sum_{m, s\in\Z}
z^mq^{m^2/4-Lm}\frac{1}{(q)_\infty}
\qbin{2L}{s}_{q^{-1}}q^{ms}.
\en
Changing $s$ to $s+L-j/2$, and 
changing $m$ to $m-2s+j$ afterwards, we obtain 
\be
\frac{1}{(q)_\infty}\sum_{m\in\Z}z^mq^{m^2/4}
\sum_{s\in\Z}\qbin{\scriptstyle2L}{\scriptstyle s+L-\frac j2}_{q^{-1}}
z^{-2s+j}q^{-(s-j/2)^2}.
\en
Changing $q$ to $q^{-1}$ we obtain the desired formula.
\end{proof}
Letting $L\rightarrow \infty$ with $L\in\Z$, 
we obtain from \eqref{twochirality} 
a formal but suggestive expression
\be
\chi_i(q^{-1},z)\chi_0(q,z)
\en
for the would-be `character' of $\Zhc^{(i)}$. 

We refer the reader to \cite{Ber} for similar results
on product formulas for shifted characters.


\section{Form factors and $\infty$-cycles}\label{formfactor}
In this section we discuss the relation
between the space of $\infty$-cycles and 
form factors of the $SU(2)$ invariant Thirring model 
(ITM).

\subsection{Form factor axioms} 

First we recall the setting. 
Let $V=\mathbb{C}v_{+} \oplus \mathbb{C}v_{-}$ be 
the vector representation of 
$\mathfrak{sl}_{2}=\C E\oplus\C F \oplus \C H$. 
The action is given by 
\be 
Ev_{+}=0, \, Ev_{-}=v_{+}, \, Fv_{+}=v_{-}, \, Fv_{-}=0, \, Hv_{\pm}=\pm v_{\pm}. 
\en 
Denote by $P \in {\rm End}(V^{\otimes 2})$ 
the permutation operator 
$P(u \otimes v)=v \otimes u$. 
The S-matrix of ITM is 
the linear operator acting on $V^{\otimes 2}$ 
defined by 
\be 
S(\beta)=
\frac{\zeta(-\beta)}{\zeta(\beta)}
\frac{\beta-\pi i P}{\beta -\pi i}. 
\en 
Here $\zeta(\beta)$ is a certain meromorphic scalar
function that accounts for the overall normalization. 
The precise formula can be found e.g. in \cite{NT}, eq.(16). 

With each local field $\mathcal{O}$ in the theory 
is associated a tower of functions 
${\bf f}^{\mathcal{O}}=
(f_{n}(\beta_{1}, \ldots , \beta_{n}))_{n \ge 0}$ 
called the form factor of $\mathcal{O}$. 
The function $f_{n}(\beta_{1}, \ldots , \beta_{n})$ 
takes values in the tensor product $V^{\otimes n}$. 
Form factor $(f_{n})_{n \ge 0}$ 
should satisfy the following axioms: 
\be
{\rm Axiom~1:} &&\!\!\!
f_{n}(\ldots , \beta_{j+1}, \beta_{j}, \ldots)=
P_{j, j+1}S_{j, j+1}(\beta_{j}-\beta_{j+1})
f_{n}(\ldots , \beta_{j}, \beta_{j+1}, \ldots), \\
{\rm Axiom~2:} &&\!\!\!
f_{n}(\beta_{1}, \ldots , \beta_{n-1}, \beta_{n}+2\pi i) \\
&& {}=
e^{\frac{n\pi i}{2}}P_{n, n-1} \cdots P_{2, 1}
f_{n}(\beta_{n}, \beta_{1}, \ldots , \beta_{n-1}), \\
{\rm Axiom~3:} &&\!\!\!
{\rm res}_{\beta_{n}=\beta_{n-1}+\pi i}f_{n}(\beta_{1}, \ldots , \beta_{n})
\\
&& {}=
(I+e^{-\frac{n\pi i}{2}}S_{n-1, n-2}(\beta_{n-1}-\beta_{n-2})
\cdots S_{n-1, 1}(\beta_{n-1}-\beta_{1})) \\
&&\quad\times f_{n-2}(\beta_1,\ldots,\beta_{n-2})
\otimes(v_+\otimes v_--v_-\otimes v_+).
\en
Here $P_{j, j+1}$ is the permutation operator 
acting on the $j$-th and $(j+1)$-st components,
and $S_{j, j'}(\beta)$ is the operator acting 
on the tensor product of the $j$-th and $j'$-th components 
of $V^{\otimes n}$ as $S(\beta)$. 

The operators $S_{i, j}(\beta)$ and $P_{i, j}$ 
commute with the action of $\mathfrak{sl}_{2}$.  
If $(f_{n})_{n\ge 0}$ is a form factor, 
then $(Ff_{n})_{n \ge 0}$ is also a form factor. 
{}For that reason, we restrict our considerations to 
only form factors satisfying the highest weight condition 
\bea 
&& 
H f_n(\beta_{1}, \ldots , \beta_{n})=
mf_n(\beta_{1}, \ldots , \beta_{n}), 
\label{eq:weight} 
\\
&& 
E f_n(\beta_{1}, \ldots , \beta_{n})=0 
\label{eq:highest} 
\ena 
for some $m \in \mathbb{Z}_{\ge 0}$ for all $n$. 
Other form factors can be obtained from these ones 
by the action of $F \in \mathfrak{sl}_{2}$.  

\subsection{The integral formula and $\infty$-cycles} 

A large class of form factors of the 
ITM is given in terms of the hypergeometric 
integral \cite{NPT}. 
It has the following structure:
\bea 
&& 
\psi_{P}(\beta_{1}, \ldots , \beta_{n}):= 
\sum_{\# M=l}v_{M}\int_{C^{l}}\prod_{p=1}^{l}d\alpha_{p} 
\prod_{p=1}^{l}\phi(\alpha_{p}; \beta_{1}, \ldots , \beta_{n}) 
\label{def:pairing} \\ 
&&\qquad {}\times 
w_{M}(\alpha_{1}, \ldots , \alpha_{l}| \beta_{1}, \ldots , \beta_{n}) 
P(X_{1}, \ldots , X_{l}| z_{1}, \ldots , z_{n}), \nn 
\ena 
where $0\le l\le n$. 
Let us explain the notation. 

First, $v_M\in V^{\otimes n}$ denotes the vector 
\be 
v_{M}:=v_{\epsilon_{1}} \otimes \cdots \otimes v_{\epsilon_{n}} 
\in V^{\otimes n}, \quad 
\en 
where the $\epsilon_j$'s are related to the 
index set $M=\{m_{1}, \ldots , m_{l}\}$ 
($1\le m_{1}< \cdots <m_{l}\le n$) 
by $M=\{j | \epsilon_{j}=-\}$. 
Second, $\phi$ and $w_M$ are fixed functions given by
\be 
&&\phi(\alpha ; \beta_{1}, \ldots , \beta_{n}):= 
\prod_{j=1}^{n}\left( 
\frac{1}{1-e^{-(\alpha-\beta_{j})}}
\frac{\Gamma(\frac{\alpha-\beta_{j}+\pi i}{-2\pi i})}
   {\Gamma(\frac{\alpha-\beta_{j}}{-2\pi i})} 
\right), 
\\
&& 
w_{M}:={\rm Skew}_{\alpha_{1}, \ldots , \alpha_{l}}g_{M}, 
\\
&& 
g_{M}(\alpha_{1}, \ldots , \alpha_{l}| \beta_{1}, \ldots , \beta_{n}) \\ 
&& {}:= 
\prod_{p=1}^{l} \left( 
\frac{1}{\alpha_{p}-\beta_{m_{p}}} 
 \prod_{j=1}^{m_{p}-1} \frac{\alpha_{p}-\beta_{j}+\pi i}{\alpha_{p}-\beta_{j}}
\right) 
\prod_{1 \le p<p' \le l}(\alpha_{p}-\alpha_{p'}+\pi i). 
\en 
In the second line, 
${\rm Skew}_{\alpha_{1}, \ldots , \alpha_{l}}$ 
stands for the skew symmetrization with respect to 
$\alpha_{1}, \ldots , \alpha_{l}$. 
Third, 
$P \in \widehat{C}_{n, l}$ is a deformed cycle, 
wherein the variables $X_a,z_j$ are related to the 
variables $\alpha_a,\beta_j$ by 
\be 
X_{a}=e^{-\alpha_{a}}, \quad z_{j}=e^{\beta_{j}}. 
\en 
Lastly, the integration contour $C$ goes along the real axis, 
except that the simple poles of the integrand at
\be
\alpha_p=\beta_j-2\pi i \Z_{\geq0}
\en
are located below $C$, and those at
\be
\alpha_p=\beta_j-\pi i+2\pi i \Z_{\geq0}
\en
above $C$. These are the only poles of the integrand.
The integral converges absolutely if $2l \le n$. 


For a cycle $P_{n,l} \in \widehat{C}_{n, l}$ we set 
\be
f_{P_{n, l}}(\beta_{1}, \ldots , \beta_{n}):= 
c_{n, l}e^{\frac{n}{4}\sum_{j=1}^{n}\beta_{j}}
\prod_{1 \le j<j' \le n}\!\!\!\zeta(\beta_{j}-\beta_{j'}) \cdot
\psi_{P_{n, l}}(\beta_{1}, \ldots , \beta_{n}), 
\en 
where $c_{n, l}$ is a constant defined by 
\be 
&& 
c_{n, l}:=\left\{ 
\begin{array}{ll} 
\prod_{a=0}^{l}d_{2k+2a, a}, & (m=2k), \\ 
\prod_{a=1}^{l}d_{2k+2a-1, a}, & (m=2k+1), 
\end{array} 
\right. \\ 
&& 
d_{n, l}:=\frac{2\pi}{\zeta(-\pi i)}(-2\pi i)^{-l-\frac{n}{2}}l!. 
\en 

The precise connection between $\infty$-cycles and form factors
is given by the following theorem proved in \cite{NT}. 
\begin{thm} 
For an $\infty$-cycle ${\bf p}=(P_{n,l})$ of weight 
$m \in \mathbb{Z}_{\ge 0}$, 
the tower ${\bf f}_{\bf p}=(f_{P_{n, l}})$ satisfies 
the form factor axioms as well as the highest weight conditions
\eqref{eq:weight}, \eqref{eq:highest}.
\end{thm} 

A local field $\mathcal{O}$ is said to have 
Lorentz spin $s$ if 
the homogeneous property 
\be 
f_{n}^{\mathcal{O}}(\beta_{1}+\Lambda , \ldots , \beta_{n}+\Lambda)= 
e^{s\Lambda}f_{n}^{\mathcal{O}}(\beta_{1}, \ldots , \beta_{n}) 
\en 
holds for its form factor. 
It is easy to see that the 
Lorentz spin of ${\bf f}_{\bf p}$ 
is equal to $\deg{\bf p}$ 
introduced in \eqref{eq:def-degree}. 
 
Suppose ${\bf f}^{\mathcal{O}}=(f_{n})$ is a form factor 
of a local field. 
Then there exists an $N \ge 0$ such that $f_{n}=0$ if $n < N$. 
{}From Axiom 3 we have 
\be 
{\rm Axiom~3':} \quad 
{\rm res}_{\beta_{N}=\beta_{N-1}+\pi i}f_{N}(\beta_{1}, \ldots , \beta_{N})=0. 
\en 
It can be shown \cite{N} 
that the function $f_{P_{N, l}}$ 
associated with a minimal cycle $P_{N, l}$ satisfies 
${\rm \hbox{Axiom }3'}$. 
As in \cite{N,JMT}, we expect the
converse to be true, namely that 
for any $N$-minimal form factor ${\bf f}=(f_{n})_{n\ge 0}$ 
satisfying \eqref{eq:weight}, \eqref{eq:highest}, 
there exists a minimal cycle $P_{N, l} \in \widehat{W}_{N, l}$ 
such that $f_{P_{N, l}}=f_{N}$. 
Under this assumption, 
Theorem \ref{thm:main-result} states that 
any ${\bf f}$ (with fixed weight $m$) 
can be written as ${\bf f}_{\bf p}$ 
for some $\infty$-cycle 
${\bf p}\in\Zbc$. 

\subsection{Realization of the space of local operators} 

The space of $\infty$-cycles is not the same as 
the space of local fields, since there are null cycles. 
In \cite{T} it is proved that 
\bea 
(x_{0}^{-}.A_{n, l-1}+ 
(x_{0}^{-})^{(2)}.A_{n, l-2})\cap\widehat{C}_{n, l}
=\ker{\left( P_{n, l} \mapsto f_{P_{n, l}} \right)}. 
\label{eq:tarasov} 
\ena 
Now we set 
\be 
\widehat{M}_{n, l}:=\widehat{C}_{n, l}/
\left( (x_{0}^{-}.A_{n, l-1}+ 
(x_{0}^{-})^{(2)}.A_{n, l-2})\cap \widehat{C}_{n, l} \right) 
\en 
and consider the projection 
\be 
\widehat{C}_{n, l} \longrightarrow \widehat{M}_{n, l}, \quad 
P_{n, l} \mapsto \bar{P}_{n, l}. 
\en 
Set 
\be 
\Lc[m]:=\{ 
\bar{\bf p}=(\bar{P}_{n, l}) | 
{\bf p}=(P_{n, l}) \in \Zhc[m] \}. 
\en 

\begin{lem} 
If $m < 0$, then we have $\Lc[m]=\{0\}$. 
\end{lem} 
\begin{proof} 
Let ${\bf p}=(P_{n, l}) \in \Zhc[m]$ be $N$-minimal 
and set $l'=\frac{N-m}{2}$. 
Then we have $P_{N, l'} \in \widehat{W}_{N, l'}$. 
In \cite{JMT} it is proved that 
\be 
\widehat{W}_{N, l'}/\left( 
x_{0}^{-}.\widehat{W}_{N, l'-1}+ 
(x_{0}^{-})^{(2)}.\widehat{W}_{N, l'-2} \right)=0 \quad 
\hbox{if} \quad N-2l' < 0. 
\en 
Hence there exists an $\infty$-cycle
\be 
{\bf q}\in x_{0}^{-}.\Zhc + (x_{0}^{-})^{(2)}.\Zhc 
\en 
such that ${\bf p}-{\bf q}$ is $(N+2)$-minimal. 
Repeating this argument, we find that for each $n,l$ we have
$\bar{P}_{n,l}=0$. Therefore $\bar{\bf p}=0$. 
\end{proof} 

Now we set $\Lc=\oplus_{m \ge 0}\Lc[m]$. 
Denote by $\mathcal{F}$ the space of form factors 
satisfying \eqref{eq:weight} and \eqref{eq:highest}. 
Then the following map is well defined: 
\be 
\boldsymbol\Phi : 
\Lc \longrightarrow \mathcal{F}, \quad 
\bar{\bf p}=(\bar{P}_{n, l}) \mapsto 
{\bf f}_{\bf p}=(f_{P_{n, l}}). 
\en 
We conjecture that the map ${\bf \Phi}$ is an isomorphism.  

In the end, let us give some examples
discussed in \cite{NT}. 
\medskip

\noindent{\bf Identity operator}.\quad 
The form factor of the identity operator $\mathbb{I}$ is 
of weight zero and zero minimal. 
It is given by 
\be 
{\bf f}^{\mathbb{I}}=(1, 0, 0, \ldots ). 
\en 
Note that $f_{2}=0$ and $f_{0} \not=0$
does not violate Axiom 3, since 
the latter becomes trivial for $n=2$, 
\be 
{\rm res}_{\beta_{2}=\beta_{1}+\pi i}f_{2}(\beta_{1}, \beta_{2})= 
(I+(-1)I)(f_{0} \otimes 
(v_{+}\otimes v_{-}-v_{-}\otimes v_{+}))=0. 
\en 
The form factor ${\bf f}^{\mathbb{I}}$ is 
obtained from the $\infty$-cycle 
\be 
{\bf 1}_{0}=(1_{0}, X, X \wedge X^{3}, 
X \wedge X^{3} \wedge X^{5}, \ldots ). 
\en 
This can be seen from \eqref{eq:tarasov} and 
the following lemma. 
\begin{lem}\label{lem:onetime-integration} 
Set 
\be 
\bar{A}_{2l , l}=\{ 
P(X_{1}, \ldots , X_{l}) \in A_{2l, l} | 
P(X_{1}, \ldots , X_{l-1}, 0)=0\}. 
\en 
Then 
\be 
\bar{A}_{2l, l} \subset x_{0}^{-}. A_{2l, l-1}+ 
(x_{0}^{-})^{(2)}. A_{2l, l-2}. 
\en 
\end{lem} 
For the proof, see Appendix \ref{app:2}.

\medskip

\noindent{\bf $su(2)$ currents}.\quad 
The $\infty$-cycles associated with $su(2)$ 
currents $j_{\sigma}^{+}$ ($\sigma=\pm$)
are of weight two and two minimal. 
They are given by 
\be 
&& 
(j_{+}^{+})_{2l+2, l}=(-1)^{l}(\prod_{j=1}^{2l+2}z_{j}^{-1}) 
X \wedge X^{3} \wedge \cdots \wedge X^{2l-1}, \\ 
&& 
(j_{-}^{+})_{2l+2, l}=X^{3} \wedge X^{5} \wedge \cdots \wedge X^{2l+1} \quad (l \ge 0). 
\en 
It is easy to see that 
\be 
j_{+}^{+}=-x_{-1}^{+}\mathbb{I}, \quad 
j_{-}^{+}=x_{1}^{+}\mathbb{I}. 
\en 

\medskip

\noindent{\bf Energy-momentum tensor}.\quad 
Denote by $T_{z}$ and $T_{\bar{z}}$ the holomorphic and 
anti-holomorphic part of the energy momentum tensor, 
respectively. 
They are of weight zero and two minimal. 
The form factors are obtained 
from the following sequences of deformed cycles: 
\be 
&& 
(T_{z})_{2l, l}=(\sum_{j=1}^{2l}z_{j})
X^{0} \wedge X^{3} \wedge X^{5} \wedge \ldots \wedge X^{2l-1}, \\ 
&& 
(T_{\bar{z}})_{2l ,l}=(-1)^{l-1}(\prod_{j=1}^{2l}z_{j}^{-1})
(\sum_{j=1}^{2l}z_j^{-1}) 
X^{0} \wedge X^{1} \wedge X^{3} \wedge \cdots \wedge X^{2l-3} \quad (l \ge 1). 
\en 
Note that these sequences are 
not $\infty$-cycles. 
However the following holds modulo 
null cycles \eqref{eq:tarasov}:
\be 
T_{z}=-ix_{1}^{-}x_{1}^{+}\mathbb{I}, \quad 
T_{\bar{z}}=-ix_{-1}^{+}x_{-1}^{-}\mathbb{I}. 
\en 
These equalities can be checked 
by using Lemma \ref{lem:onetime-integration}.

\appendix 

\section{Polynomial realization of evaluation modules}
\label{app:1}

In this appendix, we derive an action of the quantum algebra 
$U_{\sqrt{-1}}$ on the space $A_n=\oplus_{l=0}^nA_{n,l}$.
The definition of the action is given in subsection \ref{action-on-cycle}.
Here we give the 
origin of the definition and some details of its derivation.

Recall that the space $A_{n,l}$ is the space of skew-symmetric polynomials in
the variables $X_1,\ldots,X_l$ of degree less than or equal to $n-1$
with coefficients in $K_n=\C(z_1,\ldots,z_n)$. We use the wedge notation
defined in \eqref{eq:def-wedge} for skew-symmetric polynomials.
For $1\leq a\leq n$ we set
\[
G_a(X)=\prod_{j=1}^{a-1}(1+z_jX)\prod_{j=a+1}^n(1-z_jX).
\]
First note the following simple fact
(see Remark below eq.(3.12) in \cite{JMT}). 
\begin{lem}\label{L1}
The skew-symmetric polynomials
$G_{p_1}\wedge\cdots\wedge G_{p_l}\,(1\leq p_1<\cdots<p_l\leq n)$
constitute a $K_n$-basis of $A_{n,l}$.
\end{lem}

Set $V=\C v_+\oplus\C v_-$. Define operators $\sigma^z,\sigma^\pm$ acting
on $V$:
\[
\sigma^z=\left(\begin{matrix}1&0\\0&-1\end{matrix}\right),
\sigma^+=\left(\begin{matrix}0&1\\0&0\end{matrix}\right),
\sigma^-=\left(\begin{matrix}0&0\\1&0\end{matrix}\right).
\]
When we consider the tensor product of $V$,
we denote by $\sigma^*_a$ the operator $\sigma^*$ acting on the $a$-th
tensor component. 
Extending the coefficient ring we set
$\widetilde{K}_n=K_n\otimes\C[q,q^{-1}]$, 
$\widetilde{A}_n=A_n\otimes\C[q,q^{-1}]$. 
We define the action of $g\in\Ures$ (see Section 3.1)
on $V^{\otimes n}\otimes \widetilde{K}_n$ by
$(\rho_{z_1}\otimes\cdots\otimes\rho_{z_n})\circ\Delta^{(n-1)}(g)$
where $\rho_z$ is the evaluation 
representation of $\Ures$ on $V\otimes\C[q,q^{-1}][z,z^{-1}]$ 
such that
\begin{equation}\label{eq:eval}
e_1\mapsto\sigma^+,\,f_1\mapsto\sigma^-,\,t_1\mapsto q^{\sigma^z},
e_0\mapsto z\sigma^-,\,f_0\mapsto z^{-1}\sigma^+,t_0\mapsto q^{-\sigma^z} 
\end{equation}
and $\Delta$ is the coproduct defined by 
\be 
\Delta(e_{i})=e_{i}\otimes 1 +t_{i} \otimes e_{i}, \quad 
\Delta(f_{i})=f_{i} \otimes t_{i}^{-1}+1\otimes f_{i}, \quad 
\Delta(t_{i})=t_{i} \otimes t_{i}. 
\en 
Here we use the Chevalley generators $e_i,f_i,t_i$, i.e.,
\[
x^+_0=e_1,x^+_{-1}=t_1^{-1}f_0,x^-_0=f_1,x^{-}_1=e_0t_1,t_0=t_1^{-1}.
\]
We induce an action of $\Ures$ on $\widetilde{A}_n$ from the action on 
$V^{\otimes n}\otimes \widetilde{K}_n$ by using an isomorphism
between $V^{\otimes n}\otimes \widetilde{K}_n$ and 
$\widetilde{A}_n$. The isomorphism is given
as follows. Let $\psi_1,\ldots,\psi_n$ be a set of Grassmann variables.
We denote by $\Lambda_n$ the Grassmann algebra generated by them.
It is an irreducible module over the fermion algebra $\Psi_n$ generated by
$\psi_a,\psi^*_a\,(1\leq a\leq n)$ such that $[\psi_a,\psi^*_b]_+=\delta_{a,b},
[\psi_a,\psi_b]_+=[\psi^*_a,\psi^*_b]_+=0$.

There are isomorphisms of vector spaces over $\widetilde{K}_n$:
\begin{equation}\label{ISOM}
V^{\otimes n}\otimes \widetilde{K}_n\simeq\Lambda_n\otimes 
\widetilde{K}_n\simeq \widetilde{A}_n.
\end{equation}
The first isomorphism is given by the Jordan-Wigner transformation
\be 
\psi^*_a=\sigma^+_a\prod_{j={a+1}}^n(-i\sigma^z_j),
\quad 
\psi_a=\sigma^-_a\prod_{j={a+1}}^n(i\sigma^z_j) 
\en 
and the identification
of $v_+^{\otimes n}\in V^{\otimes n}$ with $1\in\Lambda_n$. The second
isomorphism is given by the identification of the left multiplication
of $\psi_a$ on $\Lambda_n\otimes \widetilde{K}_n$ with the wedge product
$G_a\wedge$ on $\widetilde{A}_n$ (see Lemma \ref{L1}).

Our goal is to compute the actions of the operators
$\X^-_{>0}(t)$, $\X^-_{\leq0}(t)$, $\X^+_{\geq0}(t)$, $\X^+_{<0}(t)$,
$\X^-_{>0}(t)^{(2)}$, $\X^-_{\leq0}(t)^{(2)}$, $\X^+_{\geq0}(t)^{(2)}$,
$\X^+_{<0}(t)^{(2)}$, $a_\pm(t)$ (see Section \ref{action-on-cycle})
in the limit $\varepsilon\rightarrow0$ by setting $q=ie^\varepsilon$.
The computation is elementary but long. We do not give all the details of
computation.

By a straightforward calculation we obtain
\begin{lem}
We have the expansion of the end term of the half current $\X^-_{>0}(t)$
up to the order $\varepsilon$.
\begin{equation}
q^{-1}x^-_1=\sum_az_a\psi_a+\left(\sum_{a<b}z_a\psi_a\sigma^z_b\right)\varepsilon
+O(\varepsilon^2).\label{E1}
\end{equation}
Similarly, we have the expansion of the end term of the half current $a_+(t)$
up to the order $\varepsilon^2$.
\begin{equation}
iq^{-1}a_1=\alpha_0+\alpha_1\varepsilon+\alpha_2\varepsilon^2
+O(\varepsilon^3),\label{E2}
\end{equation}
where
\be
&& 
\alpha_0=\sum_az_a, \quad \alpha_1=-\sum_az_a\sigma^z_a+4\sum_{a<b}z_a\psi_a\psi^*_b, \\ 
&& 
\alpha_2-\alpha_0/2=4\sum_{a<b<c}z_a\psi_a\sigma^z_b\psi^*_c.
\en
\end{lem}

The other terms in the half current $\X^-_{>0}(t)$ are determined recursively
by the equation
\begin{equation}
\X^-_{>0}(t)-q^{-1}tx^-_1=\frac{it}{[2]}[iq^{-1}a_1,\X^-_{>0}(t)].
\label{E5}
\end{equation}
Using (\ref{E1}) and (\ref{E2}), one can check
\begin{prop}\label{P1}
We have the expansion
\begin{equation}
\X^-_{>0}(t)=\gamma_0(t)+\gamma_1(t)\varepsilon+O(\varepsilon^2),
\end{equation}
where
\begin{equation}
\gamma_0(t)=\sum_aA_a(t)\psi_a,\gamma_1(t)=
\sum_{a<b}B_{a,b}(t)\psi_a\sigma^z_b+\sum_{a<b<c}C_{a,b,c}(t)\psi_a\psi_b\psi^*_c,
\end{equation}
and
\begin{eqnarray*}
&&A_a(t)=\frac{z_at}{1-z_at}\prod_{j>a}\frac{1+z_jt}{1-z_jt},\quad
B_{a,b}(t)=\frac{1-z_bt}{1+z_bt}A_a(t),\\
&&C_{a,b,c}(t)=8C_{a,b}(t)\frac{A_c(t)}{z_ct},\quad
C_{a,b}(t)=\frac{z_at}{1-z_at}\left(\prod_{a<j<b}\frac{1+z_jt}{1-z_jt}\right)
\frac{z_bt}{1-z_bt}.
\end{eqnarray*}
\end{prop}
The following equality is useful in this calculation.
\begin{equation}
\sum_{b>a}A_b(t)=\frac12\left(\prod_{j>a}\frac{1+z_jt}{1-z_jt}-1\right).
\end{equation}

Instead of repeating similar calculations for $\X^-_{\leq0}(t)$,
$\X^+_{\geq0}(t)$, $\X^+_{<0}(t)$, we can use symmetries. Let
$\alpha$ be an anti-algebra map of the algebra
$\Psi_n\otimes \widetilde{K}_n$ given by
\begin{equation}
\alpha(z_a)=z_{n+1-a},\quad\alpha(\psi_a)=z_{n+1-a}^{-1}\psi^*_{n+1-a},\quad
\alpha(\psi^*_a)=z_{n+1-a}\psi_{n+1-a}.
\end{equation}
We extend it to the space $\Psi_n\otimes \widetilde{K}_n[[t, t^{-1}]]$
of formal series in $t$ by setting $\alpha(t)=t$.
Similarly, we define an algebra map $\beta$ by
\begin{equation}
\beta(z_a)=z_{n+1-a}^{-1},\quad\beta(\psi_a)=(-1)^{a+1}\psi^*_{n+1-a},\quad
\beta(\psi^*_a)=(-1)^{a+1}\psi_{n+1-a},
\end{equation}
and extend it to formal series in $t$ by setting $\beta(t)=t^{-1}$.
Then we have $\alpha^2={\rm id}$, $\beta^2=(-1)^{(n+1)\epsilon}{\rm id}$,
$\alpha\beta=(-1)^{(n+1)\epsilon}\beta\alpha$, where $\epsilon$ is the parity
in the fermion algebra $\Psi_n$. Namely, $\epsilon=0$ on the even part of
$\Psi_n$, and $\epsilon=1$ on the odd part. We also have $\alpha(T)=T$, 
$\beta(T)=T^{-1}$, where 
\be 
T=\prod_{a=1}^{n}(i\sigma^z_a). 
\en 

\begin{prop}\label{P2}
We have
\begin{eqnarray}
\X^+_{\geq0}(t)&=&-iTt^{-1}\alpha(\X^-_{>0}(t))\bmod\varepsilon^2,
\label{E3}\\
\X^+_{<0}(t)&=&iT^{-1}\beta(\X^-_{>0}(t))\bmod\varepsilon^2,
\\
\X^-_{\leq0}(t)&=&t\beta\alpha(\X^-_{>0}(t))\bmod\varepsilon^2.
\end{eqnarray}
\end{prop}
\begin{proof}
We use the equalities
\begin{eqnarray}
\alpha(q^{-1}a_1)=q^{-1}a_1\bmod\varepsilon^3,\label{E4}\\
\beta(q^{-1}a_1)=-qa_{-1}\bmod\varepsilon^3.
\end{eqnarray}
Let us prove (\ref{E3}). It is easy to check
\[
{}-iTt^{-1}\alpha(x^-_1q^{-1}t)=x^+_0\bmod\varepsilon^2.
\]
By using (\ref{E4}), the equation (\ref{E5}), which determines $\X^-_{>0}(t)$,
is transformed into
\[
{}-iTt^{-1}\alpha(\X^-_{>0}(t))-x^+_0=-\frac{it}{[2]}[iq^{-1}a_1,
{}-iTt^{-1}\alpha(\X^-_{>0}(t))]\bmod\varepsilon^2.
\]
{}From the defining relations of the Drinfeld currents we have
\[
\X^+_{\geq0}(t)-x^+_0=-\frac{it}{[2]}[iq^{-1}a_1,
\X^+_{\geq0}(t)].
\]
Therefore, we obtain (\ref{E3}). Other cases are similar.
\end{proof}
{}From Propositions \ref{P1} and \ref{P2}, a short calculation leads to
\begin{prop}
On the space $\Lambda_n\otimes \widetilde{K}_n$,
the actions of the half currents in the limit $q\rightarrow\sqrt{-1}$
are given as follow. We write them on the subspace isomorphic to
$A_{n,l}$ by $(\ref{ISOM})$, and show
the equalities as rational functions in $t$. However, they should be properly
understood as equalities of power series in $t$
for $\X^-_{>0}(t)$, $\X^+_{\geq0}(t)$, and in $t^{-1}$
for $\X^-_{\leq0}(t)$, $\X^+_{<0}(t)$.
\begin{eqnarray}
&&\X^-_{>0}(t)=-\X^-_{\leq0}(t)=
\sum_a \frac{z_{a}t}{1-z_{a}t}\prod_{j>a}\frac{1+z_{j}t}{1-z_{j}t} \psi_a,\\
&&\X^-_{>0}(t)^{(2)}=\X^-_{\leq0}(t)^{(2)}=
i\sum_{a<b} \frac{z_{a}t}{1-z_{a}t} 
\left( \prod_{a < j < b} \frac{1+z_{j}t}{1-z_{j}t} \right) 
\frac{z_{b}t}{1-z_{b}t} \psi_a\psi_b,
\label{E6}\\
&&\X^+_{\geq0}(t)=-\X^+_{<0}(t)=-i^{n-2l-1}\sum_a\frac1{1-z_at}\prod_{j<a}
\frac{1+z_jt}{1-z_jt}\psi^*_a,\\
&&\X^+_{\geq0}(t)^{(2)}=\X^+_{<0}(t)^{(2)}=(-1)^ni\sum_{a<b}\frac1{1-z_at}
\left(\prod_{a<j<b}\frac{1+z_jt}{1-z_jt}\right)\frac1{1-z_bt}\psi^*_a\psi^*_b.
\end{eqnarray}
\end{prop}
We set
\begin{eqnarray*}
&&b_{\pm}(t):=\sum_{\pm n>0}a_n(q^{-1}t)^n,\\
&&b^{(2)}_{\pm}(t):=\sum_{\pm n>0}\frac{a_{2n}}{q+q^{-1}}
(q^{-1}t)^{2n}.
\end{eqnarray*}
\begin{prop}
The actions of $b_\pm(t)$, $b_\pm^{(2)}(t)$ in the limit
$\varepsilon\rightarrow0$ are given as follows.
\begin{eqnarray}
b_+(t)&=&-i\sum_{m\geq1,{\rm odd}}\frac{\sum_az_a^m}mt^m,
\label{E7}\\
2b^{(2)}_+(t)&=&\hbox{the even part of }\label{E8}\\
&&\hskip-40pt\sum_a\frac{z_at}{1+z_at}\sigma^z_a
{}-4\sum_{a<b}\frac{z_at}{1+z_at}
\left(\prod_{a<j<b}\frac{1-z_jt}{1+z_jt}\right)\frac1{1+z_bt}\psi_a\psi^*_b.
\nonumber\\
&&\hskip-50ptb_-(t)=-\beta(b_+(t)),\quad
b^{(2)}_-(t)=-\beta(b^{(2)}_+(t)).\label{E9}
\end{eqnarray}
\end{prop}
\begin{proof}
We use the equality
\[
b_+(t)=\frac1{q-q^{-1}}{\rm log}(1+(q-q^{-1})t_1^{-1}
[x^+_0,\Xc^-_{>0}(t)]).
\]
A simple calculation shows (\ref{E7}) and (\ref{E8}).
For $b_-(t)$, we use the properties of the algebra map $\beta$, (\ref{E6}) and
\[
\beta(x^+_0)=-iT^{-1}x^-_0\bmod\varepsilon^2,\quad \beta(t^{-1}_1)=t_1,
\]
and obtain
\[
\beta(b_+(t))=\frac1{q-q^{-1}}
{\rm log}(1+(q-q^{-1})t_1[x^-_0,\Xc^+_{<0}(t)])
\bmod\varepsilon^2.
\]
Since the expression in the right hand side is exactly $-b_-(t)$, we have
(\ref{E9}) in the limit $\varepsilon\rightarrow0$.
\end{proof}


{\it Proof of Proposition \ref{prop:Umod}}.
The final step is to rewrite the actions by the second part of the isomorphisms
(\ref{ISOM}).
It is often useful in the work
to note that these actions are symmetric with respect to $z_1,\ldots,z_n$
on $A_n$. This fact follows from the construction in \cite{TV} as explained
in \cite{JMT}.

The actions (\ref{act1}) and (\ref{act2}) follow from the equalities
\begin{eqnarray}
&&\sum_aA_a(t)G_a(x)=F_n(t,X),\label{G1}\\
&&4\sum_{a<b}C_{ab}(t)G_a(X)\wedge G_b(X)=F_n^{(2)}(t|X_1,X_2).\label{G2}
\end{eqnarray}
Note that there are no poles at $t=X$ in $F_n(t,X)$, or $t=X_1,X_2$
in $F_n^{(2)}(t|X_1,X_2)$. By induction, we can check the equalities
at the pole $t=z_n^{-1}$ (and therefore at $t=z_a^{-1}$ for all $a$),
and at $t=\infty$.

The action (\ref{act3}) follows from (\ref{E7}). Let us prove (\ref{act4}).
We rewrite as $\sigma^z_a=1-2\psi_a\psi^*_a$ in (\ref{E8}). Then, the proof
reduces to the case $l=1$, which is equivalent to the equality
\begin{eqnarray*}
&&-\frac{z_at}{1+z_at}G_a(X)-2\sum_{b<a}\frac{z_bt}{1+z_bt}
\left(\prod_{b<j<a}\frac{1-z_jt}{1+z_jt}\right)\frac{G_b(X)}{1+z_at}\\
&&\hskip15pt=-\frac t{X+t}\left(G_a(X)-\frac{\Theta_n(X)}{\Theta_n(-t)}G_a(-t)\right).
\end{eqnarray*}

The proof of (\ref{act5}) is similar.

Let us prove (\ref{act6}). We use (\ref{E3}). note that
\[
t^{-1}\alpha(\gamma_0(t))=\sum_a\frac1{1-z_at}\prod_{j<a}\frac{1+z_at}{1-z_at}
\psi^*_a.
\]
This operator sends $G_a(X)$ to
\[
\frac1{1-z_at}\prod_{j<a}\frac{1+z_at}{1-z_at}=\frac{G_a(t)}{\Theta_n(t)}.
\]
Computing the effect of $T$, we obtain (\ref{act6}).

Similarly, the proof of (\ref{act7}) 
reduces to the equality for $b<a$, 
\be 
&& 
\sum_c{\rm res}_{u=z_c^{-1}}\frac{G_a(-u)G_b(u)-G_a(u)G_b(-u)}
{\Theta_n(-u)\Theta_n(u)}\frac{du}{u-t} \\ 
&& {}= 
\frac1{1-z_bt}\Bigl(\prod_{b<j<a}
\frac{1+z_jt}{1-z_jt}\Bigr)\frac1{1-z_at}.
\en 
The first term with $G_a(-u)G_b(u)$ vanishes because there are no poles at
$u=z_c^{-1}$ for any $c$. The second term with $-G_a(u)G_b(-u)$
can be calculated by the residue at
$u=t$ because there are no pole at $u=-z_c^{-1}$.\qed

\section{Proof of Lemma \ref{lem:onetime-integration}}
\label{app:2}

Following the line of \cite{T}, we give a proof of 
Lemma \ref{lem:onetime-integration} by using 
the fermionic realization of the space $A_{n, l}$ 
given in Appendix \ref{app:1}. 
Let $\psi_{1}, \ldots , \psi_{n}$ be a set of 
Grassmann variables. 
Denote by $K_n[\psi_{1}, \ldots , \psi_{n}]$ 
the exterior algebra generated by them 
over $K_n=\mathbb{C}(z_{1}, \ldots , z_{n})$. 
Introduce the degree defined by $\deg{\psi_{a}}=1$ and 
denote by $K_n[\psi_{1}, \ldots , \psi_{n}]_{l}$ 
the homogeneous component 
of degree $l$. 
Then we have the isomorphism 
\be 
\mathcal{C}_{n, l}: K_n[\psi_{1}, \ldots , \psi_{n}]_{l}
\buildrel\sim\over\longrightarrow 
A_{n, l}. 
\en 
See Lemma \ref{ISOM} for the definition of this isomorphism. 
The action of $x_{0}^{-}$ and $(x_{0}^{-})^{(2)}$ 
on $A_{n}=\oplus_{l=0}^{n}A_{n, l}$ is intertwined with multiplication
on $K_n[\psi_{1}, \ldots , \psi_{n}]$ with the following elements: 
\be 
x_{0}^{-} \leftrightarrow \Sigma_1=\sum_{a=1}^{n}(-1)^{n-a}\psi_{a}, \quad 
i(x_{0}^{-})^{(2)} \leftrightarrow 
\Sigma_{2}=\sum_{1 \le a<b \le n}(-1)^{a+b}\psi_{a}\psi_{b}. 
\en 

Now let us prove the lemma. 
Set $n=2l$ in the construction above. 
Introduce a set of new generators $\{\varphi_{a}\}_{a=1}^{2l}$ 
of $K_n[\psi_{1}, \ldots , \psi_{2l}]$ given by 
\be 
\varphi_{a}=\psi_{a}-\psi_{2l} \,\, (a=1, \ldots , 2l-1), \quad 
\varphi_{2l}=\psi_{2l}. 
\en 
The elements $\Sigma_{1}$ and $\Sigma_{2}$ are represented in terms of these generators 
as follows: 
\bea
&& 
\Sigma_{1}=\sum_{a=1}^{2l-1}(-1)^{a}\varphi_{a}, 
\label{eq:xzero-1} \\ 
&& 
\Sigma_{2}=\tilde{\Sigma}_{2}-\Sigma_{1}(\varphi_{2l-1}-\varphi_{2l}), 
\label{eq:xzero-2} 
\ena 
where
\be 
\tilde{\Sigma}_{2}=\sum_{1 \le a<b \le 2l-2}(-1)^{a+b}\varphi_{a}\varphi_{b}. 
\en 

{}From the definition of the isomorphism $\mathcal{C}_{2l, l}$ we can see that 
\be 
\mathcal{C}_{2l, l}(K_n[\varphi_{1}, \ldots , \varphi_{2l-1}]_{l})= 
\bar{A}_{2l,l}. 
\en 
{}From \eqref{eq:xzero-1} we have that 
\bea {}\quad
K_n[\varphi_{1}, \ldots , \varphi_{2l-1}]_{l}= 
\Sigma_{1} \cdot K_n[\varphi_{1}, \ldots , \varphi_{2l-2}]_{l-1}+ 
K_n[\varphi_{1}, \ldots , \varphi_{2l-2}]_{l}. 
\label{eq:tarasov-1} 
\ena 
In \cite{T} it is proved that the map 
\be 
K_n[\varphi_{1}, \ldots , \varphi_{2l-2}]_{l-2} \longrightarrow 
K_n[\varphi_{1}, \ldots , \varphi_{2l-2}]_{l}, \quad 
x \mapsto \tilde{\Sigma}_{2} \cdot x 
\en 
is an isomorphism. 
{}From this fact, \eqref{eq:xzero-2} and \eqref{eq:tarasov-1}, 
we find that 
\be 
K_n[\varphi_{1}, \ldots , \varphi_{2l-1}]_{l} \subset 
\Sigma_{1}\cdot K_n[\varphi_{1}, \ldots , \varphi_{2l}]_{l-1}+ 
\Sigma_{2} \cdot K_n[\varphi_{1}, \ldots , \varphi_{2l-2}]_{l-2}. 
\en 
This completes the proof. 
\qed

\bigskip 
\noindent
{\it Acknowledgments.}\quad
We thank Masaki Kashiwara for helpful discussions.
JM is partially supported by 
the Grant-in-Aid for Scientific Research (B2) no.12440039, 
and TM is partially supported by 
(A1) no.13304010, Japan Society for the Promotion of Science.
EM is partially supported by the National Science
Foundation (NSF) grant DMS-0140460.
YT is supported by the Japan Society for the Promotion of Science.

\end{document}